\documentclass[a4paper, 11pt]{article}
\usepackage[dvips]{graphicx}
\usepackage{verbatim}%comment
\usepackage{amsthm}%lemma, proposition, definition
\usepackage{amsmath,amssymb,color}
\usepackage{color} 
\RequirePackage[OT1]{fontenc}
\RequirePackage{amsthm,amsmath, amssymb, enumerate}
\RequirePackage[round, sort&compress, authoryear]{natbib}
\RequirePackage[colorlinks,citecolor=blue,urlcolor=blue]{hyperref}

\usepackage{eufrak}
\usepackage{graphicx, amsthm, url,pstricks,fancyhdr}
\RequirePackage{bbm,amsthm}
\RequirePackage{amssymb,mathrsfs}

\usepackage{float}
\RequirePackage{enumitem}
\RequirePackage{bigints}
\RequirePackage{mathtools}
\RequirePackage[stable]{footmisc}
\usepackage{bm}

\theoremstyle{plain}
\newtheorem{theorem}{Theorem}[section]                                          
                          
\newtheorem{lemma}[theorem]{Lemma}
\newtheorem{corollary}[theorem]{Corollary}

\theoremstyle{definition}

\theoremstyle{remark}
\newtheorem{remark}[theorem]{Remark}

\makeatletter \@addtoreset{equation}{section} \makeatother

\newcommand{\Prob}{\mathbb{P}\,}

\newcommand{\bea}{\begin{eqnarray}}
\newcommand{\ena}{\end{eqnarray}}
\newcommand{\beq}{\begin{equation}}
\newcommand{\enq}{\end{equation}}
\newcommand{\beas}{\begin{eqnarray*}}
\newcommand{\enas}{\end{eqnarray*}}

\newcommand{\EE}{{\mathbb{E}} }

\def\heiko{\mathrel{\raise.3ex\hbox{\scalebox{.7}{%
    \rotatebox[origin=c]{-7}{/}%
    \kern-.35em\rotatebox[origin=c]{-7}{/}}}}}%

\title{Assessing the multivariate normal approximation of the maximum likelihood estimator from high-dimensional, heterogeneous data\protect} 
\author{Andreas Anastasiou$^{1}$ 
\\
\small{$^{1}$ The London School of Economics and Political Science}
\\
}
\date{}
\begin{document}
\pagenumbering{roman}
\pagenumbering{arabic}
\maketitle

\begin{abstract}
\indent The asymptotic normality of the maximum likelihood estimator (MLE) under regularity conditions is a cornerstone of statistical theory. In this paper, we give explicit upper bounds on the distributional distance between the distribution of the MLE of a vector  parameter, and the multivariate normal distribution. We work with possibly high-dimensional, independent but not necessarily identically distributed random vectors. In addition, we obtain upper bounds in cases where the MLE cannot be expressed analytically.
\end{abstract}


%\tableofcontents

\section{Introduction}
\label{sec:intro} The assessment of the quality of various normal approximations has attracted the interest of statisticians for many years. In general this is not an easy task and as \cite{Kiefer} points out, \textit{to give explicitly useful bounds on the departure from the asymptotic normal distribution as a function of the sample size seems to be a terrifically difficult problem}. Since then, Berry-Esseen type bounds have been derived for general (mainly linear) statistics; see for example \cite{Koroljuk} for the case of \textit{U}-statistics.

Due to the fact that the Maximum Likelihood Estimator (MLE) is not in general a linear function of the random variables, it was only recently that the assessment of its asymptotic normality has started getting significant attention in statistical research. Obtaining a quantitative statement related to the normal approximation of the MLE can be helpful to assess whether using the limiting distribution is an acceptable approximation or not. In addition, such results can save both money and time by giving a good indication on whether a larger sample size is indeed necessary, for a good approximation to hold.

The case of a scalar MLE for observations from single-parameter distributions is the first that has been covered in a series of papers. The existing approaches are mainly split into two categories based on whether a powerful technique called Stein's method (as first introduced in \cite{Stein1972}) was employed in order to get distributional bounds, or not. In the former category, where Stein's method was used, one can measure the MLE-related normal approximation error in a wide range of metrics, such as Zolotarev-type distances (for example the Wasserstein distance) and the Kolmogorov metric. \cite{Anastasiou_Reinert} provide the most general approach, where bounds on the distributional distance between the distribution of the MLE and the normal distribution are given and no restrictions are imposed on the form of the MLE. \cite{Anastasiou_Ley} give a different approach to the problem based on a combination of Stein's method with the Delta method for situations where the MLE can be expressed as a function of the sum of independent terms. Their strategy consists in benefiting from this special form of the MLE, which allows the direct usage of Stein's method on a sum of random elements. The bounds given in \cite{Anastasiou_Ley} are simpler than those obtained in \cite{Anastasiou_Reinert}. We note however, that an obvious advantage of the methodology developed in \cite{Anastasiou_Reinert} is its wider applicability as it works for all MLE settings (not requiring the MLE to be of a special form) and even for cases where an analytic expression of the MLE is not known.
%In their recent paper, \cite{Anastasiou_Gaunt} use the multivariate Delta method to expand the results of \cite{Anastasiou_Ley} in the scenario of independent random vectors that follow a multi-parameter distribution.
In the recent contribution of \cite{Anastasiou_m_dependence} the independence assumption is relaxed and the normal approximation of the MLE is assessed under the presence of a local dependence structure between the random variables. The resulting Zolotarev-type bounds are of the optimal $\mathcal{O}\left(n^{-1/2}\right)$ distance, while the obtained bounds on the Kolmogorov distance are $\mathcal{O}\left(n^{-1/4}\right)$. 

In the second category, where Stein's method is not used, bounds are given in the Kolmogorov distance. Using the Delta method and under the requirement that the MLE can be expressed as a function of the sum of independent random elements, \cite{Pinelis_Molzon} provide uniform and non-uniform Berry–Esseen bounds on the rate of convergence to normality for various statistics, among which is the MLE. The conditions used are partly different than those in \cite{Anastasiou_Ley}, where the Delta method was also employed. The bounds achieve the optimal $\mathcal{O}\left(n^{-1/2}\right)$ order. \cite{Pinelis} extends the results of \cite{Pinelis_Molzon} in cases where the MLE is not necessarily a function of the sum of independent random terms. Under conditions, he shows that the MLE can be tightly enough bracketed between two smooth enough functions, which makes the Delta method applicable. With regards to the Kolmogorov distance, the obtained bounds are again of the optimal order, which is an advantage over the Stein's method related approaches of the previous paragraph, where the order of the bound on the Kolmogorov distance is only $\mathcal{O}(n^{-1/4})$. However, the results given in \cite{Anastasiou_Reinert} and in the current paper are more general in the sense that firstly, they cover a larger family of metrics (in which the bounds are of the optimal $n^{-1/2}$ order) and secondly, under assumptions, are applicable when the MLE is not known analytically.

In this paper, we give upper bounds on the distributional distance between the distribution of a vector MLE and the multivariate normal, which under specific regularity conditions (given at a later stage) is the MLE's limiting distribution. We partly employ multivariate Stein's method and our focus is on independent but not necessarily identically distributed random vectors. The bounds obtained are explicit in terms of the sample size and the parameter. We are the first to give results for situations where the vector MLE can not be expressed in a closed form. The wide applicability of the maximum likelihood estimation method adds to the importance of our results. Among others, an MLE is used in ordinary and generalised linear models, time series analysis and a large number of other situations related to hypothesis testing and confidence intervals; see Section \ref{subsec:linear_reg} for bounds related to linear regression models.

The notation which is used throughout the paper is as follows. The parameter space is $\Theta \subset \mathbb{R}^d$ equipped with the Euclidean norm. Let $\boldsymbol{\theta} = (\theta_1,\theta_2, \ldots, \theta_d)^{\intercal}$ denote a parameter from the parameter space, while $\boldsymbol{\theta_0} = \left(\theta_{0,1}, \theta_{0,2}, \ldots, \theta_{0,d}\right)^{\intercal}$ denotes the true, but unknown, value of the parameter. The probability density (or probability mass) function is denoted by $f(\boldsymbol{x}|\boldsymbol{\theta})$, where $\boldsymbol{x} = (\boldsymbol{x_1}, \boldsymbol{x_2}, \ldots, \boldsymbol{x_n}) \in \mathbb{R}^n$. The likelihood function is $L(\boldsymbol{\theta}; \boldsymbol{x}) = f(\boldsymbol{x}|\boldsymbol{\theta})$. Its natural logarithm, called the log-likelihood function is denoted by $\ell(\boldsymbol{\theta};\boldsymbol{x})$. A maximum likelihood estimate (not seen as a random vector) is a value of the parameter which maximises the likelihood function. For many models the MLE as a random vector exists and is also unique, in which case it is denoted by $\boldsymbol{\hat{\theta}_n(X)}$; see \cite{Makelainen} for a set of assumptions that ensure existence and uniqueness. This is known as the `regular' case. However, existence and uniqueness of the MLE can not be taken for granted, see e.g. \cite{Billingsley} for an example of non-uniqueness.

For $\boldsymbol{X_1}, \boldsymbol{X_2}, \ldots, \boldsymbol{X_n}$ being independent but not necessarily identically distributed (i.n.i.d.) random vectors, we denote by $f_i(\boldsymbol{x},\boldsymbol{\theta})$ the probability density (or mass) function of $\boldsymbol{X_i}$. The likelihood function is $L(\boldsymbol{\theta}; \boldsymbol{x}) = \prod_{i=1}^{n}f_i(\boldsymbol{x_i}|\boldsymbol{\theta})$. With the parameter space $\Theta$ being an open subset of $\mathbb{R}^d$, the asymptotic normality of the MLE holds under the following regularity conditions as expressed in \cite{Hoadley}:
\begin{itemize}
\item [(N1)] $\boldsymbol{\hat{\theta}_n(X)} \xrightarrow[]{{\rm p}} \boldsymbol{\theta_0}$, as $n \rightarrow \infty$, where $\boldsymbol{\theta_0}$ is the true parameter value;

\item [(N2)] the Hessian matrix $J_k(\boldsymbol{X_k}, \boldsymbol{\theta}) = \left\lbrace \frac{\partial^2}{\partial \theta_i\partial\theta_j}\log(f_k(\boldsymbol{X_k}|\boldsymbol{\theta})) \right\rbrace_{i,j = 1,2,\ldots,d} \in \mathbb{R}^{d \times d}$ and the gradient vector $\nabla (\log(f_k(\boldsymbol{X_k}|\boldsymbol{\theta})))\in \mathbb{R}^{d \times 1}$ exist almost surely $\forall k \in \left\lbrace 1,2,\ldots, n\right\rbrace$ with respect to the probability measure $\Prob$;

\item [(N3)] $J_k(\boldsymbol{X_k}, \boldsymbol{\theta})$ is a continuous function of $\boldsymbol{\theta}$, $\forall k = 1,2,\ldots,n$, almost surely with respect to $\Prob$ and is a measurable function of $\boldsymbol{X_k}$;

\item [(N4)] $\EE_{\boldsymbol{\theta}}\left[\nabla (\log(f_k(\boldsymbol{X_k}|\boldsymbol{\theta})))\right] = \boldsymbol{0}, \; k = 1,2,\ldots,n$;

\item [(N5)] with $\boldsymbol{y}^{\intercal}$ denoting the transpose of a vector $\boldsymbol{y}$, $$\EE_{\boldsymbol{\theta}}\left[\left[\nabla (\log(f_k(\boldsymbol{X_k}|\boldsymbol{\theta})))\right]\left[\nabla (\log(f_k(\boldsymbol{X_k}|\boldsymbol{\theta})))\right]^{\intercal}\right] = -\EE\left[J_k(\boldsymbol{X_k},\boldsymbol{\theta})\right] =: I_k(\boldsymbol{\theta});$$

\item [(N6)] for
\begin{equation}
\label{Fisher_i.n.i.d}
\bar{I}_n(\boldsymbol{\theta}) = \frac{1}{n}\sum_{j=1}^{n}I_j(\boldsymbol{\theta}),
\end{equation}
there exists a matrix $\bar{I}(\boldsymbol{\theta}) \in \mathbb{R}^{d \times d}$ such that $\bar{I}_n(\boldsymbol{\theta}) \xrightarrow[n \rightarrow \infty]{} \bar{I}(\boldsymbol{\theta})$. In addition, $\bar{I}_n(\boldsymbol{\theta}), \bar{I}(\boldsymbol{\theta})$ are symmetric matrices for all $\boldsymbol{\theta}$ and $\bar{I}(\boldsymbol{\theta})$ is positive definite;

\item [(N7)] for some $\delta > 0$, $\frac{\sum_{k}\EE_{\boldsymbol{\theta_0}}\left|\boldsymbol{\lambda}^{\intercal}\nabla (\log(f_k(\boldsymbol{X_k})))\right|^{2+\delta}}{n^{\frac{2+\delta}{2}}} \xrightarrow[n \rightarrow \infty]{} 0$ for all $\boldsymbol{\lambda} \in \mathbb{R}^{d}$;

\item [(N8)] with $\|.\|$ the ordinary Euclidean norm on $\mathbb{R}^d$, then for $k,i,j \in \left\lbrace 1,2,\ldots, d\right\rbrace$ there exist $\epsilon > 0$, $K>0$, $\delta > 0$ and random variables $B_{k,ij}(\boldsymbol{X_k})$ such that
\begin{itemize}
\item [(i)] $\sup\left\lbrace \left|\frac{\partial^2}{\partial \theta_i\partial\theta_j}\log(f_k(\boldsymbol{X_k}|\boldsymbol{t}))\right|:\|\boldsymbol{t} - \boldsymbol{\theta_0}\|\leq \epsilon \right\rbrace \leq B_{k,ij}(\boldsymbol{X_k})$;
\item [(ii)] $\EE\left|B_{k,ij}(\boldsymbol{X_k})\right|^{1 + \delta}\leq K$.
\end{itemize}
\end{itemize}
Assuming that $\boldsymbol{\hat{\theta}_n(X)}$ exists and is unique, the following theorem gives the result for the asymptotic normality of the MLE in the case of i.n.i.d. random vectors in a slightly different way than \cite{Hoadley}.
\vspace{0.1in}
\begin{theorem}
\label{Theorem_asymptotic_MLE_i.n.i.d}
Let $\boldsymbol{X_1}, \boldsymbol{X_2}, \ldots, \boldsymbol{X_n}$ be independent random vectors with probability density (or mass) functions $f_i(\boldsymbol{x_i}|\boldsymbol{\theta})$, where $\boldsymbol{\theta} \in \Theta \subset \mathbb{R}^d$. Assume that the MLE exists and is unique and that the regularity conditions (N1)-(N8) hold. Also let $\boldsymbol{Z} \sim {\rm N}_d\left(\boldsymbol{0},I_{d \times d}\right)$, where $\boldsymbol{0}$ is the $d \times 1$ zero vector and $I_{d \times d}$ is the $d \times d$ identity matrix. Then, for $\bar{I}_n(\boldsymbol{\theta})$ as in \eqref{Fisher_i.n.i.d}
\begin{equation}
\label{i.n.i.d_asymptotic_result}
\sqrt{n}\left[\bar{I}_n(\boldsymbol{\theta_0})\right]^{\frac{1}{2}}\left(\boldsymbol{\hat{\theta}_n(X)} - \boldsymbol{\theta_0}\right) \xrightarrow[{n \to \infty}]{{\rm d}} \boldsymbol{Z}.
\end{equation}
\end{theorem}
\begin{proof}
\cite{Hoadley} proves in Theorem 2, p.1983 that under the regularity conditions (N1)-(N8),
\begin{equation}
\nonumber \sqrt{n}\left(\boldsymbol{\hat{\theta}_n(X)} - \boldsymbol{\theta_0}\right) \xrightarrow[{n \to \infty}]{{\rm d}} \left[\bar{I}(\boldsymbol{\theta_0})\right]^{-\frac{1}{2}}\boldsymbol{Z}.
\end{equation}
Using this result and (N6) we obtain that
\begin{equation}
\nonumber \left[\bar{I}_n(\boldsymbol{\theta_0})\right]^{\frac{1}{2}}\sqrt{n}\left(\boldsymbol{\hat{\theta}_n(X)} - \boldsymbol{\theta_0}\right) \xrightarrow[{n \to \infty}]{{\rm d}} \left[\bar{I}(\boldsymbol{\theta_0})\right]^{\frac{1}{2}}\left[\bar{I}(\boldsymbol{\theta_0})\right]^{-\frac{1}{2}}\boldsymbol{Z} =  \boldsymbol{Z},
\end{equation}
which is the result of the theorem.
\end{proof}

The interest is on assessing the quality of the asymptotic normality of the MLE in \eqref{i.n.i.d_asymptotic_result}. For any three times differentiable function $h\hspace{-0.03in}:\hspace{-0.03in}\mathbb{R}^d\hspace{-0.03in}\rightarrow\hspace{-0.03in}\mathbb{R}$, we abbreviate $\|h\|:= \sup |h|, \|h\|_1 := \underset{i}{\sup}\left|\frac{\partial}{\partial x_i}h\right|, \|h\|_2 := \underset{i,j}{\sup}\left|\frac{\partial^2}{\partial x_i\partial x_j}h\right|$, and $\|h\|_3 := \underset{i,j,k}{\sup}\left|\frac{\partial^3}{\partial x_i\partial x_j\partial x_k}h\right|$. For $j \in \left\lbrace 1,2,3 \right\rbrace$, let
\begin{equation}
\label{class_multi}
\hspace{-0.02in}H \hspace{-0.03in}= \hspace{-0.03in}\left\lbrace h\hspace{-0.03in}:\hspace{-0.03in}\mathbb{R}^d \rightarrow \mathbb{R}\hspace{-0.03in}:\hspace{-0.02in}h {\rm\:is\:three\:times\:differentiable\:with\:bounded\:} \|h\|, \|h\|_j\right\rbrace
\end{equation}
be the class of test functions used in the paper. We will give upper bounds on
\begin{equation}
\label{quantityofinterest}
\left|\EE\left[h\left(\sqrt{n}\left[\bar{I}_n(\boldsymbol{\theta_0})\right]^{\frac{1}{2}}\left(\boldsymbol{\hat{\theta}_n(X)} - \boldsymbol{\theta_0}\right)\right)\right]- \EE[h(\boldsymbol{Z})]\right|,
\end{equation}
where $\boldsymbol{Z} \sim {\rm N}_d\left(\boldsymbol{0},I_{d \times d}\right)$. The bounds are explicit in terms of the sample size and $\boldsymbol{\theta_0}$. The main result of the paper is given in Theorem \ref{Theorem_i.n.i.d}, where we obtain a general upper bound on \eqref{quantityofinterest} which holds under slightly weaker assumptions than the usual, sufficient regularity conditions (N1)-(N8) used for the asymptotic normality of the MLE. The generality of the bound adds to its importance as it can be applied in various different occasions. Furthermore, Theorem \ref{THEOREMMULTIIMPLICIT} is also substantial since we achieve to obtain upper bounds related to the asymptotic normality of the MLE, even when the MLE is not known analytically, but it is assumed to be within an $\epsilon$-neighbourhood of $\boldsymbol{\theta_0}$, for $\epsilon>0$.

The paper is organised as follows. Section \ref{sec:non_identically} first treats the case of independent but not necessarily identically distributed (i.n.i.d.) random vectors. The upper bound on the distributional distance between the distribution of the vector MLE and the multivariate normal distribution is presented. Special attention is given to linear regression models with an application to the simplest case of the straight-line model, where apart from the upper bound, we also give results from a simulation study. Furthermore, under weaker regularity conditions, we explain how the bound can be simplified for the case of i.i.d. random vectors. Specific theoretical and empirical results for independent random variables from the normal distribution under canonical parametrisation are given. In Section \ref{sec:multi-parameter_implicit} we explain how the results can be expanded when no analytic expression of the vector MLE is available. We briefly illustrate the results for the Beta distribution with both shape parameters unknown. In order to make the paper easily readable, we only provide an outline of the proof of our main Theorem \ref{Theorem_i.n.i.d}, with the complete proof being given in Section \ref{sec:proofs}. In addition, some technical results and proofs of corollaries that are not essential for the smooth understanding of the paper are confined in the Appendix.
\section{Bounds for multi-parameter distributions}
\label{sec:non_identically}
In this section we examine the case of  i.n.i.d. $t$-dimensional random vectors, for $t \in \mathbb{Z}^+$. We give an upper bound on the distributional distance between the distribution of the MLE and the multivariate normal. An example from linear models then follows. The last subsection covers, under weaker regularity conditions, the case of i.i.d. random vectors and an example from the normal distribution under canonical parametrisation serves as illustration of our results. It is worth mentioning that the MLE in this example is not a sum of random variables and classical Stein method approaches cannot be applied directly.
%We have already indicated that the existence and uniqueness of the MLE is not to be taken for granted. In order to secure its existence and uniqueness when the likelihood function $L(\boldsymbol{\theta};\boldsymbol{x})$ is twice continuously differentiable varying in an open parameter space $\Theta \subset \mathbb{R}^{d}$, one can make the following assumptions from \cite{Makelainen}:
%\begin{itemize}
%\item[(A1)]
%$\underset{\boldsymbol{\theta} \to \partial \Theta}\lim L(\boldsymbol{\theta};\boldsymbol{x}) = 0$, where $\partial \Theta$ is the boundary of the parameter space,
%\item[(A2)] the Hessian matrix
%\begin{center}
%$\boldsymbol{H}(\boldsymbol{\theta};\boldsymbol{x}) = \left\lbrace\frac{\partial^2}{\partial \theta_i \partial \theta_j}L(\boldsymbol{\theta};\boldsymbol{x})\right\rbrace_{i,j=1,\dots,d}\;$ (where $\frac{\partial}{\partial \theta_i}$ denotes partial derivatives)
%\end{center}
%of second partial derivatives is negative definite at every point $\boldsymbol{\theta} \in \Theta$ for which the gradient vector
%\begin{center}
%$\nabla L(\boldsymbol{\theta};\boldsymbol{x}) = \left\lbrace \frac{\partial}{\partial \theta_i}L(\boldsymbol{\theta};\boldsymbol{x}) \right\rbrace_{i = 1,\dots,d}$
%\end{center}
%vanishes.
%\end{itemize}
\subsection{A general bound}
\label{subsec:General_inid}
The normal approximation in \eqref{i.n.i.d_asymptotic_result} is an asymptotic result and our motivation is to assess the quality of this normal approximation through explicit, for finite sample size, upper bounds on the distributional distance of interest. From now on, $\bar{I}_n(\boldsymbol{\theta})$ is as in \eqref{Fisher_i.n.i.d}. Let the subscript $(m)$ denote an index for which the quantity $\left|\hat{\theta}_n(\boldsymbol{x})_{(m)} - \theta_{0,(m)}\right|$ is the largest among the $d$ components;
\begin{align}
%\label{(m)}
\nonumber & (m) \in \left\lbrace 1,\ldots,d\right\rbrace\;{\rm is\;such\;that\;} \left|\hat{\theta}_n(\boldsymbol{x})_{(m)} - \theta_{0,(m)}\right| \geq \left|\hat{\theta}_n(\boldsymbol{x})_j - \theta_{0,j}\right|, \forall j \in \left\lbrace 1,\ldots, d \right\rbrace.
\end{align}
For ease of presentation, let us introduce the following notation:
\begin{equation}
\begin{aligned}
\label{cm}
& Q_{(m)}=Q_{(m)}(\boldsymbol{X},\boldsymbol{\theta_0}) := \hat{\theta}_n(\boldsymbol{X})_{(m)} - \theta_{0,(m)}\\
& Q_{j}=Q_{j}(\boldsymbol{X},\boldsymbol{\theta_0}) := \hat{\theta}_n(\boldsymbol{X})_{j} - \theta_{0,j}, \quad \forall j \in \left\lbrace 1,2,\ldots,d \right\rbrace\\
& T_{lj} = T_{lj}\left(\boldsymbol{\theta_0},\boldsymbol{X}\right) = \frac{\partial^2}{\partial\theta_l\partial\theta_j}\ell(\boldsymbol{\theta_0};\boldsymbol{X}) + n[\bar{I}_n(\boldsymbol{\theta_0})]_{lj}, \quad j,l \in \left\lbrace 1,2,\ldots,d \right\rbrace\\
& \tilde{V} = \tilde{V}(n,\boldsymbol{\theta_0}) := \left[\bar{I}_n(\boldsymbol{\theta_0})\right]^{-\frac{1}{2}}\\
& \xi_{ij} = \frac{1}{\sqrt{n}}\sum_{k=1}^{d}\tilde{V}_{jk}\frac{\partial}{\partial \theta_k}\log(f_i(\boldsymbol{X_i}|\boldsymbol{\theta_0})), \; i \in \left\lbrace 1,2,\ldots,n \right\rbrace,\; j \in \left\lbrace 1,2,\ldots,d \right\rbrace.
\end{aligned}
\end{equation}
Notice that, using conditions (N5) and (N6), $\EE\left[T_{lj}\right] = 0$ and in general, we expect $T_{lj}$ to be small. The main result of the paper is as follows.
\vspace{0.1in}
\begin{theorem}
\label{Theorem_i.n.i.d}
Let $\boldsymbol{X_1}, \boldsymbol{X_2}, \ldots, \boldsymbol{X_n}$ be i.n.i.d. $\mathbb{R}^t$-valued, $t \in \mathbb{Z}^+$, random vectors with probability density (or mass) function $f_i(\boldsymbol{x_i}|\boldsymbol{\theta})$, for which the parameter space $\Theta$ is an open subset of $\mathbb{R}^d$. Assume that the MLE exists and is unique and that (N1)-(N6) are satisfied. In addition, assume that for any $\boldsymbol{\theta_0} \in \Theta$ there exists $0<\epsilon=\epsilon(\boldsymbol{\theta_0})$ and functions $M_{kjl}(\boldsymbol{x}),\;\forall k,j,l \in \left\lbrace 1,2,\ldots,d \right\rbrace$ such that $\left|\frac{\partial^3}{\partial \theta_{k}\partial \theta_{j}\partial \theta_{l}}\ell(\boldsymbol{\theta},\boldsymbol{x})\right| \leq M_{kjl}(\boldsymbol{x})$ for all $\boldsymbol{\theta} \in \Theta$ with $|\theta_j - \theta_{0,j}| < \epsilon$ $\forall j \in \left\lbrace 1,2,\ldots,d\right\rbrace$. Also, for $Q_{(m)}$ as in \eqref{cm}, assume that $\EE\left[\left(M_{kjv}(\boldsymbol{X})\right)^2\middle|\left|Q_{(m)}\right| < \epsilon \right] < \infty$. Let $\left\lbrace \boldsymbol{X_i'}, i=1,2,\ldots,n \right\rbrace$ be an independent copy of $\left\lbrace \boldsymbol{X_i}, i=1,2,\ldots,n\right\rbrace$. For $\boldsymbol{Z} \sim {\rm N}_{d}\left(\boldsymbol{0},I_{d\times d}\right)$, $h \in H$, where $H$ is as in \eqref{class_multi}, and with $Q_{j}$, $T_{lj}$, and $\xi_{ik}$ as in \eqref{cm}, it holds that
\begin{align}
\label{final_bound_regression}
\nonumber &\left|\EE\left[h\left(\sqrt{n}\left[\bar{I}_n(\boldsymbol{\theta_0})\right]^{\frac{1}{2}}\left(\boldsymbol{\hat{\theta}_n(X)} - \boldsymbol{\theta_0}\right)\right)\right]- \EE[h(\boldsymbol{Z})]\right|\\
& \leq \frac{1}{\sqrt{n}}\left(\|h\|_1K_1(\boldsymbol{\theta_0}) + \|h\|_2K_2(\boldsymbol{\theta_0}) + \|h\|_3K_3(\boldsymbol{\theta_0})\right) + \frac{2\|h\|}{\epsilon^2}\EE\left[\sum_{j=1}^{d}Q_j^2\right],
\end{align}
where,
\begin{align}
\nonumber & K_1(\boldsymbol{\theta_0}) = \sum_{k=1}^{d}\sum_{l=1}^{d}\left|\tilde{V}_{lk}\right|\sum_{j=1}^{d}\sqrt{\EE\left[Q_j^2\right]\EE\left[T_{kj}^2\right]}\\
\label{remainder_term_R1}
& + \frac{1}{2}\sum_{k=1}^{d}\sum_{l=1}^{d}\left|\tilde{V}_{lk}\right|\sum_{j=1}^{d}\sum_{v=1}^{d}\sqrt{\EE\left[Q_j^2Q_v^2\right]}\sqrt{\EE\left[\left(M_{kjv}(\boldsymbol{X})\right)^2\middle|\left|Q_{(m)}\right| < \epsilon \right]}\\
\label{remainder_term_R2}
& K_2(\boldsymbol{\theta_0}) = \frac{1}{4\sqrt{n}}\left\lbrace\sum_{j=1}^{d}\sqrt{\sum_{i=1}^{n}{\rm Var}\left[n\xi_{ij}^2\right]} + 2\sum_{k=1}^{d-1}\sum_{j=k+1}^{d}\sqrt{\sum_{i=1}^{n}{\rm Var}\left[n\xi_{ij}\xi_{ik}\right]}\right\rbrace\\
\label{remainder_term_R3}
& K_3(\boldsymbol{\theta_0}) = \frac{1}{12n}\sum_{i=1}^{n}\EE\left[\sum_{m=1}^{d}\left|\sum_{l=1}^{d}\tilde{V}_{ml}\left(\frac{\partial}{\partial \theta_l}\left\lbrace\log(f_i(\boldsymbol{X_i'}|\boldsymbol{\theta_0})) - \log(f_i(\boldsymbol{X_i}|\boldsymbol{\theta_0}))\right\rbrace\right)\right|\right]^3.
\end{align}
\end{theorem}
\begin{remark}
\label{remark_multi_order}
\textbf{(1)} At first glance, the bound seems complicated. However, the examples that follow show that the terms are easily calculated giving an expression for the bound, which is of the optimal $n^{-1/2}$-order.\\
\textbf{(2)} Assuming that $\bar{I}_n(\boldsymbol{\theta_0}) = \mathcal{O}(1)$ in \eqref{Fisher_i.n.i.d} yields, for fixed $d$, $\EE\left[\sum_{j=1}^{d}Q_j^2\right] = \mathcal{O}\left(n^{-1}\right)$. To see this, first use that from the asymptotic normality of the MLE as expressed in Theorem \ref{Theorem_asymptotic_MLE_i.n.i.d}, $\sqrt{n}\EE\left[\boldsymbol{\hat{\theta}_n(X)} - \boldsymbol{\theta_0}\right] \xrightarrow[{n \to \infty}]{{}} \boldsymbol{0}$ and thus 
\begin{equation}
\nonumber \EE\left[Q_j\right] = o\left(\frac{1}{\sqrt{n}}\right),\; \forall j \in \left\lbrace 1,2, \ldots, d\right\rbrace.
\end{equation}
Secondly, from Theorem \ref{Theorem_asymptotic_MLE_i.n.i.d} we also get that
\begin{equation}
\label{ordervariancei.n.i.d}
n\left[\bar{I}_n(\boldsymbol{\theta_0})\right]^{\frac{1}{2}}{\rm Cov}\left[\boldsymbol{\hat{\theta}_n(X)}\right]\left[\bar{I}_n(\boldsymbol{\theta_0})\right]^{\frac{1}{2}} \xrightarrow[{n \to \infty}]{{}} I_{d \times d}.
\end{equation}
Assuming that the matrix $\bar{I}_n(\boldsymbol{\theta_0})$ is $\mathcal{O}(1)$, it follows from \eqref{ordervariancei.n.i.d} that ${\rm Var}\left[\hat{\theta}_n(\boldsymbol{X})_j\right] = \mathcal{O}\left(n^{-1}\right),\;\forall j \in \left\lbrace 1,2,\ldots,d \right\rbrace$ and therefore,
\begin{equation}
\label{orderMSEi.n.i.d}
\EE\left[Q_j^2\right] = {\rm Var}\left[\hat{\theta}_n(\boldsymbol{X})_j\right] + \left[\EE\left[Q_j\right]\right]^2 = \mathcal{O}\left(n^{-1}\right).
\end{equation}
\textbf{(3)} With $T_{lj}$ as in \eqref{cm}, using (N5), (N6), and  the fact that $\boldsymbol{X_1}, \boldsymbol{X_2}, \ldots, \boldsymbol{X_n}$ are independent yields
\begin{align}
\label{orderTlj}
\mathbb{E}\left[T_{lj}^2\right] = \sum_{i=1}^{n}{\rm Var}\left[\frac{\partial^2}{\partial\theta_l\partial\theta_j}\log\left(f_i(\boldsymbol{X_i}|\boldsymbol{\theta_0})\right)\right],
\end{align}
%In the same way, we have that under the null
%\begin{equation}
%\label{orderTlj_null}
%\mathbb{E}\left(\left(T_{lj}^*\right)^2\right) = n{\rm Var}\left(\frac{\partial^2}{\partial\theta_{l+r}\partial\theta_{j+r}}\log\left(f(\boldsymbol{X_1}|\boldsymbol{\theta_0})\right)\right),
%\end{equation}
meaning that $\mathbb{E}\left[T_{lj}^2\right]$ is $\mathcal{O}(n)$.\\
\textbf{(4)} Using \eqref{orderMSEi.n.i.d} and \eqref{orderTlj}, then if $\bar{I}_n(\boldsymbol{\theta_0}) = \mathcal{O}(1)$ it can be deduced that
\begin{equation}
\nonumber K_1(\boldsymbol{\theta_0}) = \mathcal{O}(1),\;\; K_2(\boldsymbol{\theta_0}) = \mathcal{O}(1),\;\; K_3(\boldsymbol{\theta_0}) = \mathcal{O}(1),
\end{equation}
where $K_1(\boldsymbol{\theta_0}), K_2(\boldsymbol{\theta_0}), K_3(\boldsymbol{\theta_0})$ are as in \eqref{remainder_term_R1}, \eqref{remainder_term_R2}, \eqref{remainder_term_R3}, respectively. Hence, the upper bound in Theorem \ref{Theorem_i.n.i.d} is $\mathcal{O}\left(n^{-1/2}\right)$.\\
%textbf{(3)} Often H\"{o}lder's inequality will be used to bound the third term as the calculation of absolute third moments can be quite complicated, even for simple multi-parameter distributions.\\
\textbf{(5)} In terms of the dimensionality $d$ of the parameter, having that $\xi_{ij} = \mathcal{O}(d)$, then $K_1(\boldsymbol{\theta_0}) = \mathcal{O}\left(d^4\right), K_2(\boldsymbol{\theta_0}) = \mathcal{O}\left(d^4\right)$ and $K_3(\boldsymbol{\theta_0}) = \mathcal{O}\left(d^6\right)$ as can be deduced from \eqref{remainder_term_R1}, \eqref{remainder_term_R2} and \eqref{remainder_term_R3}, respectively. The last term of the bound in \eqref{final_bound_regression} is of order $d$ in terms of the dimensionality of the parameter. Thus, for $d\gg n$ the bound does not behave well, but $d$ could grow moderately with $n$. For example $d = o\left(n^{\alpha}\right),\;0 <\alpha < \frac{1}{12}$ would still yield a bound which goes to zero as $n$ goes to infinity.
\end{remark}
{\raggedright \textit{Outline of the proof}}.
From the definition of the MLE, $\frac{\partial}{\partial \theta_{k}} l\left(\boldsymbol{\hat{\theta}_n(x)};\boldsymbol{x}\right) = 0 \;\forall k \in \left\lbrace 1, 2, \ldots, d\right\rbrace.$ A second-order Taylor expansion of $\frac{\partial}{\partial \theta_{k}} l\left(\boldsymbol{\hat{\theta}_n(x)};\boldsymbol{x}\right)$ about $\boldsymbol{\theta_0}$ yields for $Q_j$ as in \eqref{cm}
%\begin{align}
%\nonumber 0 = &\frac{\partial}{\partial \theta_k}\ell(\boldsymbol{\theta_0;x}) + \sum_{j=1}^{d}(\hat{\theta}_n(\boldsymbol{x})_j - \theta_{0_j})\left(\frac{\partial^2}{\partial \theta_k\partial\theta_j}\ell(\boldsymbol{\theta_0};\boldsymbol{x})\right)\\
%\nonumber & + \frac{1}{2}\sum_{j=1}^{d}\sum_{q=1}^{d}(\hat{\theta}_n(\boldsymbol{x})_j - \theta_{0_j})(\hat{\theta}_n(\boldsymbol{x})_q - \theta_{0_q})\left(\frac{\partial^3}{\partial\theta_k\partial\theta_j\partial\theta_q}\ell(\boldsymbol{\theta};\boldsymbol{x})\Big|_{\substack{\boldsymbol{\theta} = \boldsymbol{\theta_0^{*}}}}\right),
%\end{align}
\begin{align}
\nonumber &\sum_{j=1}^{d}Q_j\left(\frac{\partial^2}{\partial\theta_k\partial\theta_j}\ell(\boldsymbol{\theta_0};\boldsymbol{x})\right) = -\frac{\partial}{\partial \theta_k}\ell(\boldsymbol{\theta_0;x}) - \frac{1}{2}\sum_{j=1}^{d}\sum_{q=1}^{d}Q_jQ_q\left(\frac{\partial^3}{\partial\theta_k\partial\theta_j\partial\theta_q}\ell(\boldsymbol{\theta};\boldsymbol{x})\Big|_{\substack{\boldsymbol{\theta} = \boldsymbol{\theta_0^{*}}}}\right),
\end{align}
with $\boldsymbol{\theta_0^{*}}$ between $\boldsymbol{\hat{\theta}_n(x)}$ and $\boldsymbol{\theta_0}$. Adding $\sum_{j=1}^{d}n[\bar{I}_n(\boldsymbol{\theta_0})]_{kj}Q_j$ on both sides of the above equation gives, for $T_{kj}$ as in \eqref{cm}, that
%\begin{align}
%\nonumber &\sum_{j=1}^{d}(\hat{\theta}_n(\boldsymbol{x})_j - \theta_{0_j})\left(-n[\bar{I}_n(\boldsymbol{\theta_0})]_{kj} + \frac{\partial^2}{\partial\theta_k\partial\theta_j}\ell(\boldsymbol{\theta_0};\boldsymbol{x}) + n[\bar{I}_n(\boldsymbol{\theta_0})]_{kj}\right)\\
%\nonumber & = -\frac{\partial}{\partial\theta_k}\ell(\boldsymbol{\theta_0;x}) - \frac{1}{2}\sum_{j=1}^{d}\sum_{q=1}^{d}(\hat{\theta}_n(\boldsymbol{x})_j - \theta_{0_j})(\hat{\theta}_n(\boldsymbol{x})_q - \theta_{0_q})\left(\frac{\partial^3}{\partial\theta_k\partial\theta_j\partial\theta_q}\ell(\boldsymbol{\theta};\boldsymbol{x})\Big|_{\substack{\boldsymbol{\theta} = \boldsymbol{\theta_0^{*}}}}\right)
%\end{align}
%so that
\begin{align}
\label{multiresultTaylor}
\nonumber \sum_{j=1}^{d}n[\bar{I}_n(\boldsymbol{\theta_0})]_{kj}Q_j & = \frac{\partial}{\partial\theta_k}\ell(\boldsymbol{\theta_0};\boldsymbol{x})+ \sum_{j=1}^{d}Q_jT_{kj}\\
& \qquad + \frac{1}{2}\sum_{j=1}^{d}\sum_{q=1}^{d}Q_jQ_q\left(\frac{\partial^3}{\partial\theta_k\partial\theta_j\partial\theta_q}\ell(\boldsymbol{\theta};\boldsymbol{x})\Big|_{\substack{\boldsymbol{\theta} = \boldsymbol{\theta_0^{*}}}}\right).
\end{align}
Using \eqref{multiresultTaylor}, which holds $\forall k \in \left\lbrace 1,2,\ldots,d \right\rbrace$, and with $\tilde{V}$ as in \eqref{cm},
\begin{align}
\nonumber & \sqrt{n}[\bar{I}_n(\boldsymbol{\theta_0})]^{\frac{1}{2}}(\boldsymbol{\hat{\theta}_n(x)} - \boldsymbol{\theta_0})\\
\nonumber & = \frac{\tilde{V}}{\sqrt{n}}\left\lbrace\vphantom{(\left(\sup_{\theta:|\theta-\theta_0|\leq\epsilon}\left|l^{(3)}(\theta;\boldsymbol{X})\right|\right)^2}\nabla (\ell(\boldsymbol{\theta_0};\boldsymbol{x})) +  \sum_{j=1}^{d}Q_j\left(\nabla \left(\frac{\partial}{\partial \theta_j}\ell(\boldsymbol{\theta_0};\boldsymbol{x})\right) + n[\bar{I}_n(\boldsymbol{\theta_0})]_{[j]}\right)\right.\\
\nonumber & \;\;\qquad\qquad+ \left.\frac{1}{2}\sum_{j=1}^{d}\sum_{q=1}^{d}Q_jQ_q\left(\nabla\left(\frac{\partial^2}{\partial\theta_j\partial\theta_q}\ell(\boldsymbol{\theta};\boldsymbol{x})\Big|_{\substack{\boldsymbol{\theta} = \boldsymbol{\theta_0^{*}}}}\right)\right)\vphantom{(\left(\sup_{\theta:|\theta-\theta_0|\leq\epsilon}\left|l^{(3)}(\theta;\boldsymbol{X})\right|\right)^2}\right\rbrace,
\end{align}
where $[\bar{I}_n(\boldsymbol{\theta_0})]_{[j]}$ is the $j^{th}$ column of the matrix $\bar{I}_n(\boldsymbol{\theta_0})$. The triangle inequality gives that
\begin{align}
\nonumber & \left|\EE\left[h\left(\sqrt{n}\left[\bar{I}_n(\boldsymbol{\theta_0})\right]^{\frac{1}{2}}\left(\boldsymbol{\hat{\theta}_n(X)} - \boldsymbol{\theta_0}\right)\right)\right]- \EE[h(\boldsymbol{Z})]\right|\\
\label{first_term_regression}
& \leq \left|\EE\left[h\left(\frac{\tilde{V}}{\sqrt{n}}\nabla (\ell(\boldsymbol{\theta_0};\boldsymbol{X}))\right)\right]-\EE[h(\boldsymbol{Z})]\right|\\
\label{second_term_regression}
& \;\;\;+\left|\EE\left[h\left(\sqrt{n}\left[\bar{I}_n(\boldsymbol{\theta_0})\right]^{\frac{1}{2}}\left(\boldsymbol{\hat{\theta}_n(X)} - \boldsymbol{\theta_0}\right)\right) - h\left(\frac{\tilde{V}}{\sqrt{n}}\nabla (\ell(\boldsymbol{\theta_0};\boldsymbol{X}))\right)\right]\right|.
\end{align}
Now, \eqref{first_term_regression} is based on $\nabla (\ell(\boldsymbol{\theta_0};\boldsymbol{x})) = \sum_{i=1}^{n}\nabla \left(\log(f_i(\boldsymbol{x_i}|\boldsymbol{\theta_0}))\right)$ which is a sum of independent random vectors. For this expression, a bound using Stein's method for multivariate normal approximation will be derived. In contrast, \eqref{second_term_regression} will be bounded using multivariate Taylor expansions. Technical difficulties arise as the third-order partial derivatives of the log-likelihood function may not be uniformly bounded in $\boldsymbol{\theta}$. Therefore, for $0 < \epsilon = \epsilon(\boldsymbol{\theta_0})$ we will condition on whether $\left|Q_{(m)}\right|$ as defined in \eqref{cm} is greater or less than the positive constant $\epsilon$ and each case will be treated separately by bounding conditional expectations. Known probability inequalities, such as the Cauchy-Schwarz and Markov's inequality, will be employed in order to derive the upper bounds in each case.
\subsection{Linear regression}
\label{subsec:linear_reg}
This subsection calculates the bound in \eqref{final_bound_regression} for linear regression models. The asymptotic normality of the MLE in linear regression models has been proven in \cite{Fahrmeir}. We give the example of a straight-line regression and the bound turns out to be, as expected, of order $\mathcal{O}\left(n^{-1/2}\right)$, where $n$ is the sample size. The following notation is used throughout this subsection. The vector $\boldsymbol{Y} = (Y_1, Y_2, \ldots, Y_n)^{\intercal}\in \mathbb{R}^{n\times 1}$ denotes the response variable for the linear regression, while $\boldsymbol{\beta} = (\beta_1, \beta_2, \ldots, \beta_d)^{\intercal} \in \mathbb{R}^{d\times 1}$ is the vector of the $d$ parameters and $\boldsymbol{\epsilon} = (\epsilon_1, \epsilon_2, \ldots, \epsilon_n)^{\intercal} \in \mathbb{R}^{n \times 1}$ is the vector of the error terms, which are i.i.d. random variables with $\epsilon_i \sim {\rm N}(0,\sigma^2), \forall i \in \left\lbrace 1,2,\ldots n\right\rbrace$. The true value of the unknown parameter $\boldsymbol{\beta}$ is denoted by $\boldsymbol{\beta_0} = (\beta_{0,1}, \beta_{0,2}, \ldots, \beta_{0,d})^{\intercal}\in \mathbb{R}^{d \times 1}$. The design matrix is
\begin{equation}
\nonumber X = \begin{pmatrix} 1 & x_{1,2} & \ldots & x_{1,d}\\1 & x_{2,2} & \ldots & x_{2,d}\\ \vdots & \vdots & \ddots & \vdots \\1 & x_{n,2} & \ldots & x_{n,d}
\end{pmatrix}.
\end{equation}
For the model
\begin{equation}
\nonumber \boldsymbol{Y} = X\boldsymbol{\beta} + \boldsymbol{\epsilon}
\end{equation}
the aim is to find upper bounds on the distributional distance between the distribution of the MLE, $\boldsymbol{\hat{\beta}}$, and the normal distribution. The probability density function for $Y_i$ is 
\begin{equation}
\label{density_linear_REGRESSION}
f_i(y_i|\boldsymbol{\beta}) = \frac{1}{\sqrt{2\pi\sigma^2}}\exp\left\lbrace -\frac{1}{2\sigma^2}\left(y_i - X_{[i]}\boldsymbol{\beta}\right)^2 \right\rbrace,
\end{equation}
where $X_{[i]}$ denotes the $i^{th}$ row of the design matrix.
%It has been proved in \cite{Fahrmeir} that under this linear regression setting, the conditions (N1)-(N8) are rather strong and for the case where our random variables are the response variables of generalised linear models, the following sufficient conditions have been proposed,
%\begin{itemize}
%\item [(F1)] With $\lambda_{\min}$ denoting the smallest eigenvalue of a symmetric matrix,
%\begin{equation}
%\nonumber \lambda_{\min}(I_{n}(\boldsymbol{\beta})) \rightarrow \infty
%\end{equation}
%\item [(F2)] For all $\delta > 0$,
%\begin{equation}
%\nonumber \max_{\boldsymbol{\beta}\in N_n(\delta)}\|V_n(\boldsymbol{\beta}) - I_{d \times d}\| \rightarrow 0,
%\end{equation}
%where for $n\in \left\lbrace 1,2,...\right\rbrace$, $N_n(\delta) = \left\lbrace \boldsymbol{\beta}: \|\left[I_n(\boldsymbol{\beta})\right]^{\frac{1}{2}}(\boldsymbol{\beta} - \boldsymbol{\beta_0})\| \leq \delta\right\rbrace$ and $V_n(\boldsymbol{\beta})$ is the normed information matrix.
%\end{itemize}
The parameter space $\Theta = \mathbb{R}^d$ is open and if $X^{\intercal}X$ is of full rank, the matrix $X^{\intercal}X$ is invertible and
\begin{equation}
\label{MLE_linear}
\boldsymbol{\hat{\beta}} = \left(X^{\intercal}X\right)^{-1}X^{\intercal}\boldsymbol{Y}.
\end{equation}
We now bound the corresponding distributional distance.
\vspace{0.1in}
\begin{corollary}
\label{Corollary_linear_regression}
Let $Y_i, i\in\left\lbrace 1,2,\ldots,n \right\rbrace$ be independent normal random variables with
\begin{equation}
\nonumber Y_i \sim {\rm N}\left(X_{[i]}\boldsymbol{\beta_0},\sigma^2\right),
\end{equation}
where $\sigma^2$ is known. Assume that the $d \times d$ matrix $X^{\intercal}X$ is of full rank. Let $\left\lbrace Y_i', i=1,2,\ldots,n \right\rbrace$ be an independent copy of $\left\lbrace Y_i, i=1,2,\ldots,n\right\rbrace$ and $\boldsymbol{Z} \sim {\rm N}_d(\boldsymbol{0},I_{d \times d})$ and $\bar{I}_n(\boldsymbol{\beta})$ is as in \eqref{Fisher_i.n.i.d}. Then for $h \in H$ as in \eqref{class_multi},
\begin{align}
\label{final_bound_regression_2}
\nonumber & \left|\EE\left[h\left(\sqrt{n}\left[\bar{I}_n(\boldsymbol{\beta_0})\right]^{\frac{1}{2}}\left(\boldsymbol{\hat{\beta}} - \boldsymbol{\beta_0}\right)\right)\right]- \EE[h(\boldsymbol{Z})]\right|\\
\nonumber & \leq \frac{\|h\|_2}{4}\sum_{j=1}^{d}\left[\sum_{i=1}^{n}{\rm Var}\left[\left(\sum_{k=1}^{d}\frac{X_{ik}}{\sigma}\left[\left[X^{\intercal}X\right]^{-\frac{1}{2}}\right]_{jk}\left(Y_i - \sum_{m=1}^{d}X_{im}\beta_{0,m}\right)\right)^2\right]\right]^{\frac{1}{2}}\\
\nonumber & + \frac{\|h\|_2}{2}\sum_{k=1}^{d-1}\sum_{j=k+1}^{d}\left[\vphantom{\left[\EE\left[\left(\sup_{\boldsymbol{\theta}:\left|\theta_m - \theta_{0,m}\right|\leq \epsilon}\left|\frac{\partial^3}{\partial \theta_{k}\partial \theta_{j}\partial \theta_{i}}\ell(\boldsymbol{\theta};\boldsymbol{X})\right|\right)^2\middle|\left|Q_{(m)}\right| < \epsilon \right]\right]^{\frac{1}{2}}}\sum_{i=1}^{n}{\rm Var}\left[\sum_{q=1}^{d}\sum_{v=1}^{d}\frac{X_{iq}X_{iv}}{\sigma^2}\left[\left[X^{\intercal}X\right]^{-\frac{1}{2}}\right]_{jq}\left[\left[X^{\intercal}X\right]^{-\frac{1}{2}}\right]_{kv}\left(Y_i - \sum_{m=1}^{d}X_{im}\beta_{0,m}\right)^2\right]\vphantom{\left[\EE\left[\left(\sup_{\boldsymbol{\theta}:\left|\theta_m - \theta_{0,m}\right|\leq \epsilon}\left|\frac{\partial^3}{\partial \theta_{k}\partial \theta_{j}\partial \theta_{i}}\ell(\boldsymbol{\theta};\boldsymbol{X})\right|\right)^2\middle|\left|Q_{(m)}\right| < \epsilon \right]\right]^{\frac{1}{2}}}\right]^{\frac{1}{2}}\\
&+\frac{\|h\|_3}{12}\sum_{i=1}^{n}\EE\left[\sum_{m=1}^{d}\left|\sum_{l=1}^{d}\frac{X_{il}}{\sigma}\left[\left[X^{\intercal}X\right]^{-\frac{1}{2}}\right]_{ml}\left(Y_i - Y'_i\right)\right|\right]^3.
\end{align}
\end{corollary}
%\begin{remark}
%Since the elements of $[X^{\intercal}X]^{-\frac{1}{2}}$ are $\mathcal{O}\left(\frac{1}{\sqrt{n}}\right)$, the order of the bound in \eqref{final_bound_regression_2} is also $\frac{1}{\sqrt{n}}$ as expected from Theorem \ref{Theorem_i.n.i.d}.
%\end{remark}
\begin{proof}
Using \eqref{density_linear_REGRESSION}, we can see that the Hessian matrix for the log-likelihood function does not depend on $\boldsymbol{y}$ and $\bar{I}_n(\boldsymbol{\beta_0}) = \frac{1}{n\sigma^2}X^{\intercal}X$. The result in \eqref{MLE_linear} yields
\begin{align}
\label{special_result_linear}
\nonumber \sqrt{n}\left[\bar{I}_n(\boldsymbol{\beta_0})\right]^{\frac{1}{2}}\left(\boldsymbol{\hat{\beta}} - \boldsymbol{\beta_0}\right) &= \frac{1}{\sigma}\left\lbrace \left[X^{\intercal}X\right]^{-\frac{1}{2}}X^{\intercal}\boldsymbol{Y} - \left[X^{\intercal}X\right]^{\frac{1}{2}}\boldsymbol{\beta_0} \right\rbrace\\
\nonumber & = \frac{1}{\sqrt{n}}\left[\sigma\sqrt{n}\left[X^{\intercal}X\right]^{-\frac{1}{2}}\right]\frac{1}{\sigma^2}\left(X^{\intercal}\boldsymbol{Y} -X^{\intercal}X\boldsymbol{\beta_0} \right)\\
& = \frac{1}{\sqrt{n}}\left[I_n(\boldsymbol{\beta_0})\right]^{-\frac{1}{2}}\frac{\mathrm{d}}{\mathrm{d}\boldsymbol{\beta}} \ell(\boldsymbol{\beta};\boldsymbol{y})\Big|_{\substack{\boldsymbol{\beta} = \boldsymbol{\beta_0}}}.
\end{align}
Having a closer look at the expression in \eqref{special_result_linear}, we notice that actually the quantity of interest $\left|\EE\left[h\left(\sqrt{n}[\bar{I}_n(\boldsymbol{\beta})]^{\frac{1}{2}}\left(\boldsymbol{\hat{\beta}} - \boldsymbol{\beta_0}\right)\right)\right]- \EE[h(\boldsymbol{Z})]\right|$ is equal to \eqref{first_term_regression}, with \eqref{second_term_regression} being equal to zero for this specific case of the linear regression model. Thus, using \eqref{bound_first_regression} and
\begin{equation}
\nonumber \frac{\partial}{\partial\beta_{k}}\log(f_i(Y_i|\boldsymbol{\beta_0})) = \frac{X_{ik}}{\sigma^2}\left(Y_i - \sum_{m=1}^{d}X_{im}\beta_{0,m}\right)
\end{equation}
in Theorem \ref{Theorem_i.n.i.d} yields the result of the corollary.
\end{proof}
{\raggedright{\textbf{Example: The simple linear model ($d$=2)}}}\\\\
Here, we apply the results of \eqref{final_bound_regression_2} to the case of a straight-line regression with two unknown parameters. The model is
\begin{equation}
\nonumber Y_i = \beta_1 + \beta_2(x_i - \bar{x}) + \epsilon_i ,\quad \forall i \in \left\lbrace 1,2,\ldots,n \right\rbrace.
\end{equation}
The unknown parameters $\beta_1$ and $\beta_2$ are the \textit{intercept} and \textit{slope} of the regression, respectively. As before, the i.i.d. random variables $\epsilon_i \sim {\rm N}(0,\sigma^2), \forall i \in \left\lbrace 1,2,\ldots, n\right\rbrace$. The MLE exists, it is unique and $\boldsymbol{\hat{\beta}} = \left(\bar{Y},\frac{\sum_{i=1}^{n}(x_i - \bar{x})Y_i}{\sum_{i=1}^{n}(x_i - \bar{x})^2}\right)^{\intercal}$.
\vspace{0.1in}
\begin{corollary}
\label{corollary_linear_simple}
Let $Y_1, Y_2, \ldots, Y_n$ be independent random variables with $Y_i \sim {\rm N}(\beta_1 + \beta_2(x_i - \bar{x}), \sigma^2)$. The case of $x_i = x_j, \; \forall i,j \in \left\lbrace 1,2,\ldots,n\right\rbrace$ with $i \neq j$ is excluded and for $\boldsymbol{Z}\sim {\rm N}_2(\boldsymbol{0}, I_{2 \times 2})$ and $h \in H$ as in \eqref{class_multi},
\begin{equation}
\label{bound_simple_linear}
\begin{aligned}
& \left|\EE\left[h\left(\sqrt{n}\left[\bar{I}_n(\boldsymbol{\beta_0})\right]^{\frac{1}{2}}\left(\boldsymbol{\hat{\beta}} - \boldsymbol{\beta_0}\right)\right)\right]- \EE[h(\boldsymbol{Z})]\right|\\
& \leq \frac{\|h\|_2}{4}\left(\frac{3\sqrt{2}}{\sqrt{n}} + \frac{\sqrt{2\sum_{i=1}^{n}(x_i - \bar{x})^4}}{\sum_{i=1}^{n}(x_i - \bar{x})^2}\right) + \frac{8\|h\|_3}{3\sqrt{\pi}}\left(\frac{1}{\sqrt{n}}+ \frac{\sum_{i=1}^{n}|x_i - \bar{x}|^3}{\left[\sum_{i=1}^{n}(x_i - \bar{x})^2\right]^{\frac{3}{2}}}\right).
\end{aligned}
\end{equation}
\end{corollary}
\begin{remark}
{\textbf{(1)}} The calculation of the bound is easy and relies only on simple sums. As expected, the order of the bound is $\mathcal{O}\left(n^{-1/2}\right)$, which is the optimal.\\
{\textbf{(2)}} We exclude the case of $x_i = x_j, \; \forall i,j \in \left\lbrace 1,2,\ldots,n\right\rbrace$ with $i \neq j$, only to ensure that $X^{\intercal}X$ is invertible. 
\end{remark}
\begin{proof}
We have that
\begin{equation}
\label{inverse_design}
X = \begin{pmatrix} 1 & x_1 - \bar{x}\\ 1 & x_2-\bar{x} \\ \vdots & \vdots\\1 & x_n-\bar{x}  \end{pmatrix}, \qquad X^{\intercal}X = \begin{pmatrix} n & 0\\ 0 & \sum_{i=1}^{n}(x_i-\bar{x})^2 
\end{pmatrix}.
\end{equation}
The result in \eqref{inverse_design} shows that $X^{\intercal}X$ is invertible if and only if $\sum_{i=1}^n(x_i - \bar{x})^2 \neq 0$, which holds if $x_i$'s are not all identical. The quantities of the bound in \eqref{final_bound_regression_2} are calculated for this specific case. We use that $Y_i - \beta_1 - (x_i - \bar{x})\beta_2 \stackrel{\mathclap{\normalfont\mbox{d}}}{=\joinrel=} \sigma Z_i$, where $Z_i \sim {\rm N}(0,1)$. For the first term in \eqref{final_bound_regression_2} we obtain that
\begin{align}
\label{bound_first_term_linear_2}
\nonumber& \sum_{j=1}^{2}\left[\sum_{i=1}^{n}{\rm Var}\left[\left(\sum_{k=1}^{2}\frac{X_{ik}}{\sigma}\left[\left[X^{\intercal}X\right]^{-\frac{1}{2}}\right]_{jk}\left(Y_i - \sum_{m=1}^{2}X_{im}\beta_m\right)\right)^2\right]\right]^{\frac{1}{2}}\\
\nonumber& = \sum_{j=1}^{2}\left[\sum_{i=1}^{n}{\rm Var}\left[\left(\left(\frac{X_{i1}}{\sigma}\left[\left[X^{\intercal}X\right]^{-\frac{1}{2}}\right]_{j1} + \frac{X_{i2}}{\sigma}\left[\left[X^{\intercal}X\right]^{-\frac{1}{2}}\right]_{j2}\right)\sigma Z_i\right)^2\right]\right]^{\frac{1}{2}}\\
%\nonumber & = \left[\sum_{i=1}^{n}{\rm Var}\left(\left(\frac{1}{\sqrt{n}}Z_i\right)^2\right)\right]^{\frac{1}{2}} + \left[\sum_{i=1}^{n}{\rm Var}\left(\left(\frac{x_i - \bar{x}}{\sqrt{\sum_{i=1}^{n}(x_i - \bar{x})^2}}Z_i\right)^2\right)\right]^{\frac{1}{2}}\\
\nonumber & =\frac{1}{n}\left[\sum_{i=1}^{n}{\rm Var}\left[Z_i^2\right]\right]^{\frac{1}{2}} + \frac{1}{\sum_{i=1}^{n}(x_i - \bar{x})^2}\left[\sum_{i=1}^{n}(x_i - \bar{x})^4{\rm Var}\left[Z_i^2\right]\right]^{\frac{1}{2}}\\
& = \sqrt{\frac{2}{n}} + \frac{\sqrt{2\sum_{i=1}^{n}(x_i - \bar{x})^4}}{\sum_{i=1}^{n}(x_i - \bar{x})^2}. 
\end{align}
For the second term of \eqref{final_bound_regression_2}, since $d=2$ then $k=1$, $j=2$ leading to
\begin{align}
\label{bound_second_term_linear_2}
\nonumber & \left[\sum_{i=1}^{n}{\rm Var}\left[\sum_{q=1}^{2}\sum_{v=1}^{2}\frac{X_{iq}X_{iv}}{\sigma^2}\left[\left[X^{\intercal}X\right]^{-\frac{1}{2}}\right]_{2q}\left[\left[X^{\intercal}X\right]^{-\frac{1}{2}}\right]_{1v}\left(Y_i - \sum_{m=1}^{2}X_{im}\beta_m\right)^2\right]\right]^{\frac{1}{2}}\\
\nonumber & = \left[\sum_{i=1}^{n}{\rm Var}\left[\frac{X_{i2}X_{i1}}{\sigma^2}\left[\left[X^{\intercal}X\right]^{-\frac{1}{2}}\right]_{22}\left[\left[X^{\intercal}X\right]^{-\frac{1}{2}}\right]_{11}\left(\sigma Z_i\right)^2\right]\right]^{\frac{1}{2}}\\
%\nonumber & = \left[\sum_{i=1}^{n}{\rm Var}\left(\frac{x_i - \bar{x}}{\sqrt{n\sum_{i=1}^{n}(x_i - \bar{x})^2}}Z_i^2\right)\right]^{\frac{1}{2}}= \frac{1}{\sqrt{n\sum_{i=1}^{n}(x_i - \bar{x})^2}}\left[\sum_{i=1}^{n}{\rm Var}\left(Z_i^2\right)(x_i - \bar{x})^2\right]^{\frac{1}{2}}\\
& = \frac{1}{\sqrt{n\sum_{i=1}^{n}(x_i - \bar{x})^2}}\left[\sum_{i=1}^{n}(x_i - \bar{x})^2{\rm Var}\left[Z_i^2\right]\right]^{\frac{1}{2}}= \sqrt{\frac{2}{n}}.
%\nonumber & = \left[\sum_{i=1}^{n}\frac{\left(x_i - \bar{x}\right)^2}{n\sigma^4 \sum_{i=1}^{n}(x_i - \bar{x})^2}{\rm Var}\left(\left(\sigma Z_i\right)^2\right)\right]^{\frac{1}{2}}\\
%\nonumber & = \frac{1}{\sqrt{n\sum_{i=1}^{n}(x_i - \bar{x})^2}}\left[2\sum_{i=1}^{n}(x_i - \bar{x})^2\right]^{\frac{1}{2}}
\end{align}
For the final term of \eqref{final_bound_regression_2}, because $Y_i'$ is an independent copy of $Y_i$, then \linebreak $Y_i' - Y_i \sim {\rm N}(0,2\sigma^2)$, with $\EE\left|Y_i' - Y_i\right|^3 = 8\frac{\sigma^3}{\sqrt{\pi}}$. Using that
\begin{equation}
\label{interestinginequality3}
(|a|+|b|)^3\leq 4\left(|a|^3+|b|^3\right), \;\; a,b \in \mathbb{R}
\end{equation}
yields
\begin{align}
\label{bound_third_term_linear_2}
\nonumber & \sum_{i=1}^{n}\EE\left[\sum_{m=1}^{2}\left|\sum_{l=1}^{2}\frac{X_{il}}{\sigma}\left[\left[X^{\intercal}X\right]^{-\frac{1}{2}}\right]_{ml}\left(Y_i - Y'_i\right)\right|\right]^3\\
\nonumber & = \sum_{i=1}^{n}\EE\left[\left|\left(\frac{X_{i1}}{\sigma}\left[\left[X^{\intercal}X\right]^{-\frac{1}{2}}\right]_{11}+\frac{X_{i2}}{\sigma}\left[\left[X^{\intercal}X\right]^{-\frac{1}{2}}\right]_{22}\right)\left(Y_i - Y'_i\right)\right|\right]^3\\
\nonumber & \leq \sum_{i=1}^{n}\EE\left[\left(\frac{1}{\sigma\sqrt{n}} + \frac{|x_i - \bar{x}|}{\sigma\sqrt{\sum_{i=1}^{n}(x_i - \bar{x})^2}}\right)\left|Y_i' - Y_i\right|\right]^3\\
& \leq 4\sum_{i=1}^{n}\left(\frac{8}{n^{\frac{3}{2}}\sqrt{\pi}} + \frac{8|x_i - \bar{x}|^3}{\left[\sum_{i=1}^{n}(x_i - \bar{x})^2\right]^{\frac{3}{2}}\sqrt{\pi}}\right) = \frac{32}{\sqrt{\pi}}\left(\frac{1}{\sqrt{n}}+ \frac{\sum_{i=1}^{n}|x_i - \bar{x}|^3}{\left[\sum_{i=1}^{n}(x_i - \bar{x})^2\right]^{\frac{3}{2}}}\right).
\end{align}
Summarizing, in the case of $Y_1, Y_2, \ldots, Y_n$ being independent random variables with $Y_i \sim {\rm N}\left(\beta_1 + \beta_2(x_i - \bar{x}), \sigma^2\right)$, we apply to \eqref{final_bound_regression_2} the results of \eqref{bound_first_term_linear_2}, \eqref{bound_second_term_linear_2} and \eqref{bound_third_term_linear_2} to obtain the assertion of the corollary.
\end{proof}
{\raggedright{\textbf{Empirical results}}}
\vspace{0.1in}
\\
Here, we study the accuracy of our bounds by simulations. For $n=10^j, j=3,4,5,6$, we start by generating $10^4$ trials of $n$ random independent observations, $y$, which follow ${\rm N}(\beta_1 + \beta_2(x_i - \bar{x}), \sigma^2)$, where $\beta_1 = 1, \beta_2 = 2, \sigma^2 = 1$ and each $x_i$ is sampled from the discrete uniform distribution in the set $\left\lbrace 1,2,\ldots,100 \right\rbrace$. Then $\sqrt{n}\left[\bar{I}_n(\boldsymbol{\beta_0})\right]^{\frac{1}{2}}\left(\boldsymbol{\hat{\beta}} - \boldsymbol{\beta_0}\right)$ is evaluated in each trial, which in turn gives a vector of $10^4$ values. We apply to these values the function $h(x,y) = \left(x^2+y^2+1\right)^{-1}$ and we calculate their sample mean, denoted by $\hat{\EE}\left[h\left(\sqrt{n}\left[\bar{I}_n(\boldsymbol{\beta_0})\right]^{\frac{1}{2}}\left(\boldsymbol{\hat{\beta}} - \boldsymbol{\beta_0}\right)\right)\right]$. The function $h$ is a member of the class $H$ as in \eqref{class_multi} with
\begin{equation}
\label{boundsh}
\|h\| = 1,\quad \|h\|_1 = \frac{3\sqrt{3}}{8},\quad \|h\|_2 = 2,\quad \|h\|_3 < 4.7.
\end{equation}
We use these values to calculate the bound in \eqref{bound_simple_linear}. We define $$Q_{h}(\boldsymbol{\beta_0}):=\left|\hat{\EE}\left[h\left(\sqrt{n}\left[\bar{I}_n(\boldsymbol{\beta_0})\right]^{\frac{1}{2}}\left(\boldsymbol{\hat{\beta}} - \boldsymbol{\beta_0}\right)\right)\right] - \tilde{\EE}[h(\boldsymbol{Z})]\right|,$$
where $\tilde{\EE}[h(\boldsymbol{Z})] = 0.461$ is the approximation of $\EE[h(\boldsymbol{Z})]$ up to three decimal places. We compare $Q_{h}(\boldsymbol{\beta_0})$ with the bound in \eqref{bound_simple_linear}, using the difference between their values as a measure of the error. The results are presented in Table \ref{tableresultsimplelinear} and are based on this particular function $h$, while the theoretical bounds that we have already given hold for any test function that belongs in the class $H$ defined in \eqref{class_multi}.
\begin{table}[H]
\caption{Simulation results for the simple linear model}
\vspace{0.04in}
\centering
\begin{tabular}{l|l|l|l}
	  $n$ & $Q_h(\boldsymbol{\beta_0})$ & Upper bound & Error\\
	  \hline
	  \hline
%  $10^2$ & 0.011 & 3.187 & 3.176\\
%  \hline
  $10^3$ & 0.007 & 1.002  & 0.995\\
  \hline
  $10^4$ & 0.005 & 0.319  & 0.314\\
  \hline
  $10^5$ & 0.003 & 0.101 & 0.098\\
    \hline
  $10^6$ & 0.001 & 0.032 & 0.031
  \end{tabular}
\label{tableresultsimplelinear}
  \end{table}
The table indicates that the bound and the error decrease as the sample size gets larger. When at each step we increase the sample size by a factor of ten, then the value of the upper bound drops by approximately a $\sqrt{10}$ factor, which is expected as the expression in \eqref{bound_simple_linear} is $\mathcal{O}\left(n^{-1/2}\right)$.
\subsection{Special case: Identically distributed random vectors}
\label{subsec:General_multivariate}
In this subsection we use weaker regularity conditions than (N1)-(N6) which were used in Theorem \ref{Theorem_i.n.i.d}, in order to find an upper bound in the case of independent and identically distributed random vectors. Following \cite{Davison}, we make the following assumptions:
\begin{itemize}[leftmargin=0.55in]
\item [(R.C.1)] The densities defined by any two different values of $\boldsymbol{\theta}$ are distinct;
\item [(R.C.2)] $\ell(\boldsymbol{\theta};\boldsymbol{x})$ is three times differentiable with respect to the unknown vector parameter, $\boldsymbol{\theta}$, and the third partial derivatives are continuous in $\boldsymbol{\theta}$;
\item [(R.C.3)] for any $\boldsymbol{\theta_0} \in \boldsymbol{\Theta}$ and for $\boldsymbol{\mathbb{X}}$ denoting the support of the data, there exists $\epsilon_0 > 0$ and functions $M_{rst}(\boldsymbol{x})$ (they can depend on $\boldsymbol{\theta_0}$), such that for $\boldsymbol{\theta} = (\theta_1, \theta_2, \ldots, \theta_d)$ and $r, s, t, j = 1,2,\ldots,d,$
\begin{equation}
\nonumber \left|\frac{\partial^3}{\partial \theta_r \partial \theta_s \partial \theta_t}\ell(\boldsymbol{\theta};\boldsymbol{x})\right| \leq M_{rst}(\boldsymbol{x}), \; \forall \boldsymbol{x} \in \boldsymbol{\mathbb{X}},\; \left|\theta_j - \theta_{0,j}\right| < \epsilon_0,
\end{equation}
with $\EE[M_{rst}(\boldsymbol{X})] < \infty$;
\item[(R.C.4)] for all $\boldsymbol{\theta} \in \Theta$, $\EE_{\boldsymbol{\theta}}[\ell_{\boldsymbol{X_i}}(\boldsymbol{\theta})] = 0$;
\item [(R.C.5)] the expected Fisher information matrix for one random vector $I(\boldsymbol{\theta})$ is finite, symmetric and positive definite. For $r,s=1,2,\ldots,d$, its elements satisfy
\begin{equation}
\nonumber n[I(\boldsymbol{\theta})]_{rs} = \EE\left\lbrace \frac{\partial}{\partial \theta_r} \ell(\boldsymbol{\theta};\boldsymbol{X})\frac{\partial}{\partial \theta_s}\ell(\boldsymbol{\theta};\boldsymbol{X})\right\rbrace = \EE\left\lbrace -\frac{\partial^2}{\partial \theta_r \partial \theta_s} \ell(\boldsymbol{\theta};\boldsymbol{X}) \right\rbrace.
\end{equation}
This condition implies that $nI(\boldsymbol{\theta})$ is the covariance matrix of $\nabla(\ell(\boldsymbol{\theta};\boldsymbol{x}))$.
\end{itemize}
These regularity conditions in the multi-parameter case resemble those in \cite{Anastasiou_Reinert} where the parameter is scalar. Under (R.C.1)-(R.C.5), \cite{Davison} shows that $\sqrt{n}[I(\boldsymbol{\theta_0)}]^{\frac{1}{2}}\left(\boldsymbol{\hat{\theta}_n(X)} - \boldsymbol{\theta_0}\right) \xrightarrow[{n \to \infty}]{{\rm d}} \boldsymbol{Z}.$ The upper bound on the distributional distance between the distribution of a vector MLE and the multivariate normal in the case of i.i.d. random vectors is the same as the bound in Theorem \ref{Theorem_i.n.i.d} and thus it is not given again. The bound can be simplified due to the fact that in the i.i.d. case $\bar{I}_n(\boldsymbol{\theta_0}) = I(\boldsymbol{\theta_0})$ and $f_i(\boldsymbol{x_i}) = f(\boldsymbol{x_i}), \; \forall i \in \left\lbrace 1,2, \ldots, n\right\rbrace$. In the next example of independent random variables from the normal distribution under canonical parametrisation with both natural parameters unknown, the bound can be easily calculated and it is, as expected, of the order $\mathcal{O}\left(n^{-1/2}\right)$.\\\\
{\textbf{Example: The normal distribution under canonical parametrisation}}\\\\
Many popular distributions which have the same underlying structure based on simple properties are exponential families, such as the normal, Gamma and Beta distributions; generalisations of exponential families can be found in \cite{Steffen} and \cite{Berk}. Most of the times, the interest is on working under the canonical parametrisation; the distribution of a random variable, $X$, is said to be a canonical \emph{multi-parameter exponential family distribution} if, for $\boldsymbol{\eta} \in \mathbb{R}^d$, the probability density (or mass) function is of the form
\begin{equation}
\nonumber f(x|\boldsymbol{\eta}) = {\rm exp}\left\lbrace \sum_{j=1}^{d}\eta_jT_j(x) - A(\boldsymbol{\eta}) + S(x)\right\rbrace\mathbbm{1}_{\{x \in B\}},
\end{equation}
where the set $B = \left\lbrace x:f(x|\theta)>0 \right\rbrace$ is the support of $X$ and does not depend on $\boldsymbol{\eta}$; $A(\boldsymbol{\eta})$ is a function of the parameter; $T_j(x)$ and $S(x)$ are functions only of the data. The vectors $\boldsymbol{\eta} = (\eta_1,\eta_2,\ldots,\eta_d)$ and $\boldsymbol{T(x)} = \left(T_{1}(x), T_{2}(x), \ldots, T_{d}(x)\right)$ are called the \textit{natural parameter vector} and \textit{natural sufficient statistic}, respectively. There is a number of reasons why the canonical parametrisation is more convenient. To start with, written in its canonical form, the probability density (or mass) function of an exponential family distribution has some convexity properties, which are then useful in dealing with moments and other functions of the natural sufficient statistic $\boldsymbol{T(x)}$. Furthermore, for each $j \in \left\lbrace 1,2,\ldots, d\right\rbrace$, if $X$ follows a canonical exponential family distribution, then $T_j(X)$ also follows an exponential family distribution and also
\begin{equation}
\nonumber \mathbb{E}\left[T_j(X)\right] = \frac{\partial}{\partial\eta_j}A(\boldsymbol{\eta}), \quad {\rm Cov}\left[T_{k}(X), T_{j}(X)\right] = \frac{\partial^2}{\partial\eta_k\partial\eta_j}A(\boldsymbol{\eta}), \; 1 \leq k,j \leq d.
\end{equation}
Apart from simplifying the theory and computation complexity in generalised linear models, there are other application areas, where natural exponential family distributions play a significant role. An example is the area of Gaussian graphical models (see \cite{Lauritzen_graphical} for more information) and the precision matrix estimation \citep{Massam}.

Here, we apply Theorem \ref{Theorem_i.n.i.d} in the case of $X_1, X_2, \ldots, X_n$ independent and identically distributed random variables from N$(\mu,\sigma^2)$, which is an exponential family distribution. Due to the importance, as explained above, of the natural parametrisation in exponential families, we are interested in
\begin{equation}
\label{natural_parameter}
\boldsymbol{\eta_0} = \left(\eta_1,\eta_2\right)= \left(\frac{1}{2\sigma^2},\frac{\mu}{\sigma^2}\right),
\end{equation}
which is the natural parameter vector. The MLE for $\boldsymbol{\eta_0}$ exists, it is unique and equal to $\boldsymbol{\hat{\eta}(X)} = \left(\hat{\eta}_1, \hat{\eta}_2\right)^{\intercal} = \frac{n}{\sum_{i=1}^{n}\left(X_i - \bar{X}\right)^2}\left(\frac{1}{2}, \bar{X}\right)^{\intercal}$; to see this, use the invariance property of the MLE and the result of \cite{Davison}, p.116, where the MLEs for $\mu$ and $\sigma^2$ are given. In contrast to Corollary \ref{corollary_linear_simple}, the MLE in the current example of the Gaussian distribution under canonical parametrisation is not a sum of random variables; therefore, classical Stein's method approaches, which require that the quantity of interest is a sum, cannot be employed. It appears that our results are the first that can be applied for such cases where the vector MLE has a general form in order to get upper bounds on the absolute value of the difference of expectations on the class of functions $H$ in \eqref{class_multi}. The results of \cite{Pinelis_Molzon} and \cite{Pinelis} can also be applied through the Delta method to give upper bounds only on the Kolmogorov distance though. The conditions (R.C.1)-(R.C.5) are satisfied. Corollary \ref{Corollary_multi_normal_bound} provides a bound on the distributional distance of interest and the proof is in the Appendix.
\vspace{0.1in}
\begin{corollary}
\label{Corollary_multi_normal_bound}
Let $X_1, X_2, \ldots, X_n$ be i.i.d. random variables that follow the N$(\mu, \sigma^2)$ distribution. Let $\boldsymbol{\eta_0}$ be as in \eqref{natural_parameter} and for ease of presentation, we denote $\alpha := \alpha(\eta_1,\eta_2) = \eta_1(1+\sqrt{\eta_1})^2 + \eta_2^2$ and $\beta := \beta(\eta_1,\eta_2) = \eta_1(1+\sqrt{\eta_1}) + \eta_2^2$. For $\boldsymbol{Z} \sim {\rm N}_2(\boldsymbol{0},I_{2 \times 2})$ and $h \in H$ as defined in \eqref{class_multi}, we have that for $n>9$
\begin{align}
\label{finalboundMultiNormal}
\nonumber &\left|\EE\left[h\left(\sqrt{n}[I(\boldsymbol{\eta_0})]^{\frac{1}{2}}(\boldsymbol{\hat{\eta}(X)} - \boldsymbol{\eta_0})\right)\right] - \EE[h(\boldsymbol{Z})] \right| < \frac{8\|h\|\left(\left(\eta_1^2 + \eta_2^2\right)(2n+15) + 2n\eta_1\right)}{\eta_1^2(n-3)(n-5)}\\
\nonumber & + \frac{\sqrt{2}n^{\frac{3}{2}}\|h\|_1}{\sqrt{\alpha}(n-5)(n-9)}\left\lbrace\vphantom{(\left(\sup_{\theta:|\theta-\theta_0|\leq\epsilon}\left|l^{(3)}(\theta;\boldsymbol{X})\right|\right)^2}2\sqrt{\frac{130}{\eta_1}+\frac{1473\eta_2^2}{\eta_1^2}}\left((\eta_1 + |\eta_2|)(\eta_1+3|\eta_2|+2\sqrt{\eta_1})+\eta_1\right)\right.\\
\nonumber & \left. \;\;\;\;\;\;\;\quad + \left(\frac{39\eta_2^2}{\eta_1^3}+\frac{10}{\eta_1^2}\right)\left(4\left|\eta_2^3\right| + \eta_1(2|\eta_2|+\eta_1)\left(3|\eta_2|+2+2\sqrt{\eta_1}\right)+\eta_1^{\frac{5}{2}}+\eta_1^3\right)\right.\\
\nonumber & \left. \;\;\;\;\;\;\;\quad + 156\left(\sqrt{\eta_1}+|\eta_2| + \eta_1\right)\left(1+\frac{3\left(|\eta_2|+\frac{\eta_1}{2}\right)^2}{\eta_1}\right)\vphantom{(\left(\sup_{\theta:|\theta-\theta_0|\leq\epsilon}\left|l^{(3)}(\theta;\boldsymbol{X})\right|\right)^2}\right\rbrace\\
\nonumber & + \frac{\|h\|_2}{2\sqrt{2n}\alpha}\left\lbrace\vphantom{(\left(\sup_{\theta:|\theta-\theta_0|\leq\epsilon}\left|l^{(3)}(\theta;\boldsymbol{X})\right|\right)^2}\sqrt{7}\left(\frac{\eta_2^2}{\eta_1}+1\right)\alpha + \eta_1\eta_2^2 + \frac{\beta^2}{\eta_1}\right.\\
\nonumber & \left. \;\;\;\;\;\;\;\quad + 3\sqrt{2}\eta_2\left[\left(\alpha-\eta_2^2\right)\left(5+\frac{\eta_2^2}{\eta_1}\right)^2 + \beta^2+ \left(2\sqrt{\eta_1}\eta_2+\frac{\alpha}{\eta_2}\right)^2\left(\frac{5}{3}+\frac{\eta_2^2}{\eta_1}\right)\right]^{\frac{1}{2}}\right\rbrace\\
& +  \frac{64\sqrt{2}\|h\|_3}{3\sqrt{n}\alpha^{\frac{3}{2}}}\left\lbrace\vphantom{(\left(\sup_{\theta:|\theta-\theta_0|\leq\epsilon}\left|l^{(3)}(\theta;\boldsymbol{X})\right|\right)^2}18\left(1+\frac{\eta_2^3}{2\eta_1^{\frac{3}{2}}\sqrt{\pi}}\right)\left(\eta_1^{\frac{3}{2}}\left(1+\sqrt{\eta_1}\right)^3 + \left|\eta_2\right|^3\right) + \frac{\eta_1^3|\eta_2|^3 + \beta^3}{\sqrt{\pi}\eta_1^{\frac{3}{2}}}\vphantom{(\left(\sup_{\theta:|\theta-\theta_0|\leq\epsilon}\left|l^{(3)}(\theta;\boldsymbol{X})\right|\right)^2}\right\rbrace.
\end{align}
\end{corollary}
\vspace{0.1in}
\begin{remark}
\textbf{(1)} The rate of convergence of the upper bound in \eqref{finalboundMultiNormal} is $\frac{1}{\sqrt{n}}$. Although the bound might seem complicated, the proof of the corollary shows that what is required for the derivation of the bound is basic calculation of expectations.\\
\textbf{(2)} As already mentioned, this example consists an indication of the advantages of our method in comparison to classical multivariate Stein's method results, which require that the quantity of interest is a sum of random variables. This is not the case in Corollary \ref{Corollary_multi_normal_bound} because $\hat{\eta}_1 = \frac{n}{2\sum_{i=1}^{n}\left(X_i - \bar{X}\right)^2}$.
\end{remark}
\vspace{0.1in}
{\raggedright{\textbf{Empirical results}}}
\\
We carry out a large-scale simulation study to investigate the accuracy of the bound in \eqref{finalboundMultiNormal}. The procedure is similar to the one followed previously when we obtained empirical results related to the example of the simple linear model in Corollary \ref{corollary_linear_simple}. Therefore, we start by generating $10^4$ trials of $n$ random independent observations, $y$, following ${\rm N}(\mu, \sigma^2)$, and the vector parameter of interest is $\boldsymbol{\eta_0} = (\eta_1,\eta_2)$ as in \eqref{natural_parameter}. We take $\mu =1, \sigma^2 = 1$ for our simulations. Then $\sqrt{n}\left[I(\boldsymbol{\eta_0})\right]^{\frac{1}{2}}\left(\boldsymbol{\hat{\eta}(X)} - \boldsymbol{\eta_0}\right)$ is evaluated in each trial, which in turn gives a vector of $10^4$ values. The function $h(x,y) = \left(x^2+y^2+1\right)^{-1}$, which belongs in the class $H$ as in \eqref{class_multi}, is then applied to these values in order to get the sample mean, denoted by $\hat{\EE}\left[h\left(\sqrt{n}\left[I(\boldsymbol{\eta_0})\right]^{\frac{1}{2}}\left(\boldsymbol{\hat{\eta}(X)} - \boldsymbol{\eta_0}\right)\right)\right]$. Using \eqref{boundsh}, we calculate the bound in \eqref{finalboundMultiNormal}. We define $$Q_{h}(\boldsymbol{\eta_0}):=\left|\hat{\EE}\left[h\left(\sqrt{n}\left[I(\boldsymbol{\eta_0})\right]^{\frac{1}{2}}\left(\boldsymbol{\hat{\eta}} - \boldsymbol{\eta_0}\right)\right)\right] - \tilde{\EE}[h(\boldsymbol{Z})]\right|,$$
where $\tilde{\EE}[h(\boldsymbol{Z})] = 0.461$ is the approximation of $\EE[h(\boldsymbol{Z})]$ up to three decimal places. We compare $Q_{h}(\boldsymbol{\eta_0})$ with the bound in \eqref{finalboundMultiNormal}, using the difference between their values as a measure of the error. The results from the simulations are shown in Table \ref{tableresultsmultinormal} below.
\begin{table}[H]
\caption{Simulation results for the ${\rm N}(1,1)$ distribution under a canonical parametrisation}
\vspace{0.04in}
\centering
\begin{tabular}{l|l|l|l}
	  $n$ & $Q_h(\boldsymbol{\eta_0})$ & Upper bound & Error\\
	  \hline
	  \hline
 % $10^3$ & 0.011 & 94.731 & 94.720\\
 % \hline
  $10^4$ & 0.010 & 29.898  & 29.888\\
  \hline
  $10^5$ & 0.009 & 9.452  & 9.443\\
  \hline
  $10^6$ & 0.006 & 2.988 & 2.982
  \end{tabular}
\label{tableresultsmultinormal}
  \end{table}
As in the results of Table \ref{tableresultsimplelinear}, we also see here that the bound and the error decrease as the sample size gets larger. To be more precise, when at each step we increase the sample size by a factor of ten, the value of the upper bound drops by a factor close to $\sqrt{10}$, which is expected since the order of the bound is $\mathcal{O}\left(n^{-1/2}\right)$, as can be seen from \eqref{finalboundMultiNormal}. In this example, the bounds are not as small as in Table \ref{tableresultsimplelinear}, with the reason being that the expression for the bound in \eqref{finalboundMultiNormal} is the result of a series of simplifications in order to obtain a relatively compact representation; see the proof of Corollary \ref{Corollary_multi_normal_bound} in the Appendix for the exact steps that lead to the expression in \eqref{finalboundMultiNormal}. The bound has more of a conceptual character and better constants are possible at the cost of a more complicated expression.
\section{Bounds when the MLE is not known explicitly}
\label{sec:multi-parameter_implicit}
The general bound of Theorem \ref{Theorem_i.n.i.d} includes terms of which the calculation requires to know an analytic expression for the MLE. This fact generates problems in models where no closed-form solution to the maximization problem is known or available; in these cases, a numerical method, such as the Newton-Raphson algorithm, can often be used to approximate the MLE and a normal approximation is still of interest. In this section, we will first explain how, under some further assumptions, we can put the dependence of the bound on the MLE only through the MSE, $\EE\left[\sum_{j=1}^{d}Q_j^2\right]$. Then, the MSE will get bounded by a quantity which is independent of $\boldsymbol{\hat{\theta}_n(X)}$ and it can therefore be used to get upper bounds on the distributional distance of interest that can be applied when the vector MLE is not expressed in a closed-form. To the best of our knowledge, such bounds have not appeared before in the literature for the case of a vector MLE that can not be expressed in a closed form. The extra assumptions are
\begin{itemize}[leftmargin=0.53in]
\item [(Con.1)] For an $\epsilon_0 = \epsilon_0(\boldsymbol{\theta_0}) > 0$, the MLE is within an $\epsilon$- neighbourhood of $\boldsymbol{\theta_0}$, in the sense that $\forall j \in \left\lbrace 1,2,\ldots,d\right\rbrace$, $\left|\hat{\theta}_n(\boldsymbol{X})_j - \theta_{0,j}\right| < \epsilon_0$;
\item [(Con.2)] for all $\boldsymbol{\theta_0} \in \Theta$, where $\Theta$ is the open parameter space,
\begin{equation}
\nonumber \sup_{\substack{\boldsymbol{\theta}:\left|\theta_q - \theta_{0,q}\right| < \epsilon_0\\ \forall q \in \left\lbrace 1,2,\ldots, d\right\rbrace}}\left| \frac{\partial{^3}}{\partial \theta_k\partial \theta_j\partial \theta_i}\log f(\boldsymbol{x_1}|\boldsymbol{\theta})\right| \leq M_{kji},
\end{equation}
where $M_{kji} = M_{kji}(\boldsymbol{\theta_0})$ is a constant that may depend only on $\boldsymbol{\theta_0}$;
\item [(Con.3)] the Hessian matrix of the second-order partial derivatives of the log-likelihood function is symmetric and invertible.
\end{itemize}
Section \ref{sec:non_identically} gave an upper bound for the distributional distance between the distribution of the MLE and the multivariate normal distribution. As explained in the outline of the proof of Theorem \ref{Theorem_i.n.i.d}, this bound in \eqref{final_bound_regression} can be split into terms coming from Stein's method, and terms due to Taylor expansions and conditional expectations. With $\tilde{V}$ as in \eqref{cm}, for ease of presentation we abbreviate the terms coming from Stein's method by 
\begin{align}
\label{D}
\nonumber & D = D(\boldsymbol{\theta_0},h,\boldsymbol{X}) := \frac{\|h\|_2}{4\sqrt{n}}\sum_{j=1}^{d}\left[{\rm Var}\left[\left(\sum_{k=1}^{d}\tilde{V}_{jk}\frac{\partial}{\partial\theta_{k}}\log f(\boldsymbol{X_1}|\boldsymbol{\theta_0})\right)^2\right]\right]^{\frac{1}{2}}\\
\nonumber& + \frac{\|h\|_2}{2\sqrt{n}}\sum_{k=1}^{d-1}\sum_{j=k+1}^{d}\left[\vphantom{(\left(\sup_{\theta:|\theta-\theta_0|\leq\epsilon}\left|l^{(3)}(\theta;\boldsymbol{X})\right|\right)^2}{\rm Var}\left[\vphantom{(\left(\sup_{\theta:|\theta-\theta_0|\leq\epsilon}\left|l^{(3)}(\theta;\boldsymbol{X})\right|\right)^2}\sum_{q=1}^{d}\sum_{v=1}^{d}\tilde{V}_{jq}\frac{\partial}{\partial\theta_{q}}\log f(\boldsymbol{X_1}|\boldsymbol{\theta_0})\tilde{V}_{kv}\frac{\partial}{\partial\theta_{v}}\log f(\boldsymbol{X_1}|\boldsymbol{\theta_0})\vphantom{(\left(\sup_{\theta:|\theta-\theta_0|\leq\epsilon}\left|l^{(3)}(\theta;\boldsymbol{X})\right|\right)^2}\right]\vphantom{(\left(\sup_{\theta:|\theta-\theta_0|\leq\epsilon}\left|l^{(3)}(\theta;\boldsymbol{X})\right|\right)^2}\right]^{\frac{1}{2}}\\
& + \frac{\|h\|_3}{12\sqrt{n}}\EE\left[\sum_{i=1}^{d}\left|\sum_{l=1}^{d}\tilde{V}_{il}\left(\frac{\partial}{\partial \theta_l}\log f(\boldsymbol{X_1'}|\boldsymbol{\theta_0}) - \frac{\partial}{\partial \theta_l}\log f(\boldsymbol{X_1}|\boldsymbol{\theta_0})\right)\right|\right]^3.
\end{align}
We will now first explain how we can put the dependence of the general bound in \eqref{final_bound_regression} on MLE only through the quantity $\EE\left[\sum_{j=1}^dQ_j^2\right]$ with $Q_j$ as in \eqref{cm}. After that, we will give an upper bound for $\EE\left[\sum_{j=1}^dQ_j^2\right]$.
\vspace{0.1in}
\\
\textbf{A bound depending on the mean squared error:} Under (Con.1) and with $\tilde{V}$, $Q_{(m)}$, $Q_j$ and $T_{kj}$ as in \eqref{cm}, and for $D$ in \eqref{D}, we obtain, using \eqref{final_bound_regression}, that 
\begin{align}
\nonumber & \left|\EE\left[h\left(\sqrt{n}[I(\boldsymbol{\theta_0})]^{\frac{1}{2}}(\boldsymbol{\hat{\theta}_n(X)} - \boldsymbol{\theta_0})\right)\right] - \EE[h(\boldsymbol{Z})] \right|\leq D\\
\label{boundDimplicit1}
&+ \frac{\|h\|_1}{\sqrt{n}}\sum_{k=1}^{d}\sum_{l=1}^{d}\left|\tilde{V}_{lk}\right|\sum_{j=1}^{d}\left[\EE\left[Q_j^2\right]\EE\left[T_{kj}^2\right]\right]^{\frac{1}{2}}\\
\label{boundDimplicit2}
& +\frac{\|h\|_1}{2\sqrt{n}}\left\lbrace\sum_{k=1}^{d}\sum_{l=1}^{d}\left|\tilde{V}_{lk}\right|\EE\left|\sum_{j=1}^{d}\sum_{i=1}^{d}Q_jQ_i\frac{\partial^3}{\partial \theta_k \partial\theta_j \partial\theta_i}\ell(\boldsymbol{\theta_0^{*}};\boldsymbol{X})\right|\right\rbrace.
\end{align}
{\textbf{Step 1: Upper bound for}} \eqref{boundDimplicit1}. Since $\EE\left[T_{kj}\right] = 0$, $\forall j, k \in \left\lbrace 1,2,\ldots, d \right\rbrace$,
\begin{align}
\label{first_term_implicit_multiparameter}
\nonumber & \eqref{boundDimplicit1} = \|h\|_1\sum_{k=1}^{d}\sum_{l=1}^{d}\left|\tilde{V}_{lk}\right|\sum_{j=1}^{d}\sqrt{\EE\left[Q_j^2\right]}\sqrt{{\rm Var}\left[\frac{\partial^2}{\partial \theta_k\partial \theta_j}\log f(\boldsymbol{X_1}|\boldsymbol{\theta_0})\right]}\\
& \leq \|h\|_1\sum_{k=1}^{d}\sum_{l=1}^{d}\left|\tilde{V}_{lk}\right|\sum_{j=1}^{d}\sqrt{\EE\left[Q_j^2\right]}\sqrt{\sum_{i=1}^{d}{\rm Var}\left[\frac{\partial^2}{\partial \theta_k\partial \theta_i}\log f(\boldsymbol{X_1}|\boldsymbol{\theta_0})\right]},
\end{align}
where the inequality comes from the trivial bound $${\rm Var}\left[\frac{\partial^2}{\partial \theta_k\partial \theta_j}\log f(\boldsymbol{X_1}|\boldsymbol{\theta_0})\right] \leq \sum_{i=1}^{d}{\rm Var}\left[\frac{\partial^2}{\partial \theta_k \partial \theta_i}\log f(\boldsymbol{X_1}|\boldsymbol{\theta_0})\right]$$ since the variance of a random variable is always non-negative. Now, using that \linebreak $\left(\sum_{j}^{d}\alpha_j\right)^2 \leq d\left(\sum_{j=1}^{d}\alpha_j^2\right)$ for $\alpha_j \in \mathbb{R}$, yields
\begin{equation}
\nonumber \left(\sum_{j=1}^{d}\sqrt{\EE\left[Q_j^2\right]}\right)^2\leq d\sum_{j=1}^{d}\EE\left[Q_j^2\right].
\end{equation} 
Taking square roots in both sides of the above inequality and applying this result to \eqref{first_term_implicit_multiparameter} gives
\begin{align}
\label{boundDimplicit1bound}
& \eqref{boundDimplicit1} \leq \|h\|_1\sqrt{d}\sum_{k=1}^{d}\sum_{l=1}^{d}\left|\tilde{V}_{lk}\right|\sqrt{\sum_{i=1}^{d}{\rm Var}\left[\frac{\partial^2}{\partial \theta_k\partial \theta_i}\log f(\boldsymbol{X_1}|\boldsymbol{\theta_0})\right]}\sqrt{\EE\left[\sum_{j=1}^{d}Q_j^2\right]}.
\end{align}
{\textbf{Step 2: Upper bound for} \eqref{boundDimplicit2}}. Notice that from (Con.2), $\left|\frac{\partial^3}{\partial \theta_k \partial\theta_j \partial\theta_i}\ell(\boldsymbol{\theta_0^{*}};\boldsymbol{x})\right| = \left|\sum_{l=1}^{n}\frac{\partial^3}{\partial \theta_k \partial\theta_j \partial\theta_i}\log f(\boldsymbol{x_l}|\boldsymbol{\theta_0^{*}})\right| \leq nM_{kji}$. Also,
\begin{align}
\nonumber &\sum_{j=1}^{d}\sum_{i=1}^{d}\left|Q_jQ_i\right|M_{kji} = \sum_{j=1}^{d}Q_j^2M_{kjj} + 2\sum_{i=1}^{d-1}\sum_{j=i+1}^{d}\left|Q_j\right|\left|Q_i\right|M_{kij}.
\end{align}
Using now that $2\alpha\beta \leq \alpha^2 + \beta^2, \forall \alpha, \beta \in \mathbb{R}$,
\begin{align}
\label{mid_step_boundDimplicit2}
\nonumber \sum_{j=1}^{d}\sum_{i=1}^{d}\left|Q_jQ_i\right|M_{kji} & \leq \sum_{j=1}^{d}Q_j^2M_{kjj} + \sum_{i=1}^{d-1}\sum_{j=i+1}^{d}\left[Q_j^2 + Q_i^2\right]M_{kji} = \sum_{j=1}^{d}Q_j^2\sum_{i=1}^{d}M_{kji}\\
& \leq \sum_{j=1}^{d}Q_j^2\sum_{m=1}^{d}\sum_{i=1}^{d}M_{kmi}.
\end{align}
Using \eqref{mid_step_boundDimplicit2} yields
\begin{align}
\label{boundDimplicit2bound}
\eqref{boundDimplicit2}
%\leq \frac{\|h\|_1}{2\sqrt{n}}\sum_{k=1}^{d}\sum_{l=1}^{d}\left|\tilde{V}_{lk}\right|\EE\left(n\sum_{j=1}^{d}\sum_{i=1}^{d}M_{kij}\left|\hat{\theta}_n(\boldsymbol{X})_j - \theta_{0,j}\right|\left|\hat{\theta}_n(\boldsymbol{X})_i - \theta_{0,i}\right|\middle|\left|Q_{(m)}\right| < \epsilon \right)\\
& \leq \frac{\|h\|_1\sqrt{n}}{2}\sum_{k=1}^{d}\sum_{l=1}^{d}\left|\tilde{V}_{lk}\right|\sum_{m=1}^{d}\sum_{i=1}^{d}M_{kmi}\EE\left[\sum_{j=1}^{d}Q_j^2\right].
\end{align}
Hence, from \eqref{boundDimplicit1bound} and \eqref{boundDimplicit2bound},
\begin{align}
\label{final_bound_implicit}
\nonumber & \left|\EE\left[h\left(\sqrt{n}[I(\boldsymbol{\theta_0})]^{\frac{1}{2}}(\boldsymbol{\hat{\theta}_n(X)} - \boldsymbol{\theta_0})\right)\right] - \EE[h(\boldsymbol{Z})] \right|\leq D\\
\nonumber &\; + \|h\|_1\sqrt{d}\sum_{k=1}^{d}\sum_{l=1}^{d}\left|\tilde{V}_{lk}\right|\sqrt{\sum_{i=1}^{d}{\rm Var}\left[\frac{\partial^2}{\partial \theta_k\partial \theta_i}\log f(\boldsymbol{X_1}|\boldsymbol{\theta_0})\right]}\sqrt{\EE\left[\sum_{j=1}^{d}Q_j^2\right]}\\
& \; + \frac{\|h\|_1\sqrt{n}}{2}\sum_{k=1}^{d}\sum_{l=1}^{d}\left|\tilde{V}_{lk}\right|\sum_{m=1}^{d}\sum_{i=1}^{d}M_{kmi}\EE\left[\sum_{j=1}^{d}Q_j^2\right].
\end{align}
Since $D$ as defined in \eqref{D}, is not related to the MLE, the upper bound in \eqref{final_bound_implicit} depends on $\boldsymbol{\hat{\theta}_n(X)}$ only through $\EE\left[\sum_{j=1}^{d}Q_j^2\right]$. Our purpose now is to find a bound for $\EE\left[\sum_{j=1}^{d}Q_j^2\right]$ that does not contain any terms related to $\boldsymbol{\hat{\theta}_n(X)}$. 
\vspace{0.1in}
\\
{\textbf{A bound on the mean squared error:}} In order to give an upper bound when $\boldsymbol{\hat{\theta}_n(X)}$ is not known explicitly but (Con.1)-(Con.3) are satisfied, we bound $\EE\left[\sum_{j=1}^{d}Q_j^2\right]$, for $Q_j$ as in \eqref{cm}, by a quantity which does not require knowledge of the MLE. The result is given in Theorem \ref{THEOREMMULTIIMPLICIT} below, followed by the proof.
\vspace{0.1in}
\begin{theorem}
\label{THEOREMMULTIIMPLICIT}
Let $\boldsymbol{X_1}, \boldsymbol{X_2}, \ldots, \boldsymbol{X_n}$ be i.i.d. $\mathbb{R}^t$-valued random elements, for  $t \in \mathbb{N}$, with probability density (or mass) function $f(\boldsymbol{x_i}|\boldsymbol{\theta})$, where $\boldsymbol{\theta}$ is the d-valued vector parameter. Assume that (Con.1), (Con.3) are satisfied. We assume existence and uniqueness of $\boldsymbol{\hat{\theta}_n(X)}$. For $J(\boldsymbol{x},\boldsymbol{\theta}) = \left\lbrace \frac{\partial^2}{\partial\theta_i\partial\theta_j}\ell\left(\boldsymbol{\theta};\boldsymbol{x}\right)\right\rbrace_{i,j=1,2,\ldots,d}$, the Hessian matrix, it holds that
\begin{align}
\label{multiMSEBOUND}
\EE\left[\sum_{j=1}^{d}Q_j^2\right] &\leq \EE\left[\sum_{k=1}^d\sum_{q=1}^d\frac{\partial}{\partial\theta_k}\ell(\boldsymbol{\theta_0};\boldsymbol{X})\frac{\partial}{\partial\theta_q}\ell(\boldsymbol{\theta_0};\boldsymbol{X}) \mathop{\sup_{\boldsymbol{\theta}:\left|\theta_j - \theta_{0,j}\right|<\epsilon}}_{\forall j \in \left\lbrace 1,2,\ldots, d\right\rbrace}\left\lbrace\left[J^{-2}(\boldsymbol{X},\boldsymbol{\theta})\right]_{kq}\right\rbrace \right] := U_1.
\end{align}
\end{theorem}
\vspace{0.1in}
\begin{proof}
From the definition of the MLE, we have that $\frac{\partial}{\partial\theta_k}\ell\left(\boldsymbol{\hat{\theta}_n(x)};\boldsymbol{x}\right) = 0$, $\forall k \in \left\lbrace1,2,\ldots,d\right\rbrace$. A first-order Taylor expansion of $\frac{\partial}{\partial\theta_k}\ell\left(\boldsymbol{\hat{\theta}_n(x)};\boldsymbol{x}\right)$ about $\boldsymbol{\theta_0}$ leads to
\begin{equation}
\label{boundMSEstart}
\sum_{j=1}^d\left(\hat{\theta}_n(\boldsymbol{x})_j - \theta_{0,j}\right)\frac{\partial^2}{\partial\theta_k\theta_j}\ell(\boldsymbol{\tilde{\theta}};\boldsymbol{x}) = -\frac{\partial}{\partial\theta_k}\ell(\boldsymbol{\theta_0};\boldsymbol{x}),
\end{equation}
where $\boldsymbol{\tilde{\theta}}$ is between $\boldsymbol{\theta_0}$ and $\boldsymbol{\hat{\theta}_n(x)}$. Since the result in \eqref{boundMSEstart} holds $\forall k \in \left\lbrace 1,2,\ldots, d\right\rbrace$, we deduce that
\begin{equation}
\nonumber \boldsymbol{\hat{\theta}_n(x)} - \boldsymbol{\theta_0} = -\left[J(\boldsymbol{\tilde{\theta}};\boldsymbol{x})\right]^{-1}\nabla\left(\ell(\boldsymbol{\theta_0};\boldsymbol{x})\right)
\end{equation}
and therefore
\begin{equation}
\nonumber \left(\boldsymbol{\hat{\theta}_n(x)} - \boldsymbol{\theta_0}\right)^{\intercal}\left(\boldsymbol{\hat{\theta}_n(x)} - \boldsymbol{\theta_0}\right) = \left[\nabla\left(\ell(\boldsymbol{\theta_0};\boldsymbol{x})\right)\right]^{\intercal}\left[J(\boldsymbol{\tilde{\theta}};\boldsymbol{x})\right]^{-2}\left(\ell(\boldsymbol{\theta_0};\boldsymbol{x})\right).
\end{equation}
Going a step further and using (Con.1), we get that
\begin{align}
\nonumber & \EE\left[\left(\boldsymbol{\hat{\theta}_n(X)} - \boldsymbol{\theta_0}\right)^{\intercal}\left(\boldsymbol{\hat{\theta}_n(X)} - \boldsymbol{\theta_0}\right)\right]\\
\nonumber & \leq \EE\left[\left[\nabla\left(\ell(\boldsymbol{\theta_0};\boldsymbol{X})\right)\right]^{\intercal}\mathop{\sup_{\boldsymbol{\theta}:\left|\theta_j - \theta_{0,j}\right|<\epsilon}}_{\forall j \in \left\lbrace 1,2,\ldots, d\right\rbrace}\left\lbrace\left[J(\boldsymbol{\theta};\boldsymbol{X})\right]^{-2}\right\rbrace\left(\ell(\boldsymbol{\theta_0};\boldsymbol{X})\right)\right],
\end{align}
which finishes the proof.
\end{proof}
\begin{remark}
\label{remark_implicit}
\textbf{(1)}  As the bound \eqref{multiMSEBOUND} does not include $\boldsymbol{\hat{\theta}_n(X)}$, in cases where a closed-form expression for the vector MLE is not available, we can still get an upper bound on the distributional distance between the distribution of the MLE and the $d$-variate standard normal, under the assumptions (R.C.1)-(R.C.5) and (Con.1)-(Con.3). Combining the results in \eqref{final_bound_implicit} and \eqref{multiMSEBOUND} and for $D$ as in \eqref{D} and $U_1$ as in \eqref{multiMSEBOUND}, we obtain that
\begin{align}
\label{final_upper_multi}
\nonumber &\left|\EE\left[h\left(\sqrt{n}[I(\boldsymbol{\theta_0})]^{\frac{1}{2}}(\boldsymbol{\hat{\theta}_n(X)} - \boldsymbol{\theta_0})\right)\right] - \EE[h(\boldsymbol{Z})] \right|\leq D\\
\nonumber & \;\;+ \|h\|_1\sqrt{dU_1}\sum_{k=1}^{d}\sum_{l=1}^{d}\left|\tilde{V}_{lk}\right|\sqrt{\sum_{i=1}^{d}{\rm Var}\left[\frac{\partial^2}{\partial \theta_k\partial \theta_i}\log f(\boldsymbol{X_1}|\boldsymbol{\theta_0})\right]}\\
& \;\; + \frac{\|h\|_1\sqrt{n}}{2}U_1\sum_{k=1}^{d}\sum_{l=1}^{d}\left|\tilde{V}_{lk}\right|\sum_{m=1}^{d}\sum_{i=1}^{d}M_{kmi}.
\end{align}
\textbf{(2)} In the special case where the second-order partial derivatives of the log-likelihood function do not depend on $\boldsymbol{x}$, then the result can be simplified, since in such scenarios $J(\boldsymbol{\theta};\boldsymbol{X}) = -n[I(\boldsymbol{\theta})]$ and $\left[J(\boldsymbol{\theta};\boldsymbol{X})\right]^{-2} = \frac{1}{n^2}[I(\boldsymbol{\theta})]^{-2}$. Applying this on \eqref{multiMSEBOUND}, leads to
\begin{equation}
\label{special_case}
U_1 = \frac{1}{n^2}\sum_{k=1}^d\sum_{q=1}^d\mathop{\sup_{\boldsymbol{\theta}:\left|\theta_j - \theta_{0,j}\right|<\epsilon}}_{\forall j \in \left\lbrace 1,2,\ldots, d\right\rbrace}\left\lbrace\left[I^{-2}(\boldsymbol{\theta})\right]_{kq}\right\rbrace\EE\left[\frac{\partial}{\partial\theta_k}\ell(\boldsymbol{\theta_0};\boldsymbol{X})\frac{\partial}{\partial\theta_q}\ell(\boldsymbol{\theta_0};\boldsymbol{X})\right].
\end{equation}
\end{remark}
\vspace{0.1in}
{\raggedright{\textbf{Example: The Beta distribution}}}\\\\
Here, we briefly explain how we can calculate $U_1$ in \eqref{multiMSEBOUND} for the specific example of i.i.d. random variables from the Beta distribution with both shape parameters unknown. An analytic expression for the MLE is not available. Let $\Psi_j(.)$ to be the $j^{th}$ derivative of the digamma function $\Psi$, with $\Psi(z) = \frac{\Gamma'(z)}{\Gamma(z)}, z>0$. The function $\Psi_j(z)$ can be defined through a sum, with
\begin{equation}
\label{psi_m}
\Psi_m(z) = (-1)^{m+1}m!\sum_{k=0}^{\infty}\frac{1}{(z+k)^{m+1}}, {\;\rm for\;} z \in \mathbb{C}\setminus\{\mathbb{Z}^{-}_{0}\}{\rm\;and\;m>0}.
\end{equation}
Corollary \ref{Cor_multi_Beta} gives the bound $U_1$ for the MSE in the case of the Beta distribution. The proof is given in the Appendix. For ease of presentation, and for $x,y>0$, allow us from now on to denote by 
\begin{align}
\label{delta}
\nonumber & \delta_I := \delta_I(\alpha,\beta) = \Psi_1(\alpha)\Psi_1(\beta) - \Psi_1(\alpha + \beta)\left(\Psi_1(\alpha) + \Psi_1(\beta)\right)\\
& C_1(x,y) := \Psi_1(x) - \Psi_1(x + y).
\end{align}
\begin{corollary}
\label{Cor_multi_Beta}
Let $X_1, X_2, \ldots, X_n$ be i.i.d. random variables from the Beta$(\alpha,\beta)$ distribution with $\boldsymbol{\theta_0} = (\alpha, \beta)$. Under (Con.1)-(Con.3) and with $U_1$ as in \eqref{multiMSEBOUND} and $\delta_I, C_1(x,y)$ as in \eqref{delta}, we get that
\begin{align}
\label{Beta_bound_MSE}
\nonumber U_1 = & \frac{1}{n[\delta_I(\alpha+\epsilon,\beta+\epsilon)]^2}\left\lbrace\vphantom{(\left(\sup_{\theta:|\theta-\theta_0|\leq\epsilon}\left|l^{(3)}(\theta;\boldsymbol{X})\right|\right)^2}C_1(\alpha,\beta)\left[(\alpha+\epsilon)^2\Psi_2^2(\beta - \epsilon) + \Psi_1^2(\alpha + \beta - 2\epsilon)\right]\right.\\
\nonumber & \left. + C_1(\beta,\alpha)\left[(\beta+\epsilon)^2\Psi_2^2(\alpha - \epsilon) + \Psi_1^2(\alpha + \beta - 2\epsilon)\right]\right.\\
& \left. + 2\Psi_1(\alpha+\beta)\left[\left(\beta+\epsilon\right)\Psi_2(\alpha-\epsilon) + \left(\alpha+\epsilon\right)\Psi_2(\beta-\epsilon)\right]\vphantom{(\left(\sup_{\theta:|\theta-\theta_0|\leq\epsilon}\left|l^{(3)}(\theta;\boldsymbol{X})\right|\right)^2}\right\rbrace .
\end{align}
\end{corollary}
%The following lemma gives a useful approach to finding a square root of a matrix. The proof of the lemma is given in \cite{Somayya}.
%\vspace{0.05in}
%\begin{lemma}
%\label{lemma_square_root}
%Let $M =\begin{pmatrix}
%A & B\\ C & D
%\end{pmatrix}$, where $A, B, C, D$ may be real or complex numbers. In addition, let $\tau$ and $\delta$ be the trace and determinant of $M$, respectively. Let $s$ and $t$, such that $s^2=\delta$ and $t^2 = \tau + 2s$. Then, if $t\neq 0$, a square root of $M$ is
%\begin{equation}
%\nonumber R = \frac{1}{t}\begin{pmatrix}
%A + s & B\\C & D+s
%\end{pmatrix}.
%\end{equation}
%\end{lemma}
%\begin{remark}
%The lemma holds for $s = \pm \sqrt{\delta}$ and $t = \pm \sqrt{\tau + 2s}$. Using different signs for those quantities gives different roots of the matrix $M$. From now on, unless otherwise stated, we take the positive roots for $s$ and $t$.
%\end{remark}
\vspace{0.1in}
\begin{remark}
This bound basically relies on the calculation of the expressions defined in \eqref{multiMSEBOUND}. It can be easily seen that it is of order $\mathcal{O}\left(n^{-1}\right)$. We deduce that, if we use the result in \eqref{Beta_bound_MSE} in order to calculate the bound in \eqref{final_upper_multi} for the specific case of the Beta distribution, then the obtained bound will be, as expected, of order $\mathcal{O}\left(n^{-1/2}\right)$.
%textbf{(3)} Often H\"{o}lder's inequality will be used to bound the third term as the calculation of absolute third moments can be quite complicated, even for simple multi-parameter distributions.\\
\end{remark}
\section{Proof of Theorem \ref{Theorem_i.n.i.d}}
\label{sec:proofs}
In this section, the complete steps of the proof of the main theorem of our paper are given. The following lemma (special case of Chebyshev's `other' inequality) is useful for bounding conditional expectations, which sometimes can be difficult to derive. The proof is given in the Appendix.
\vspace{0.1in}
\begin{lemma}
\label{Lemmaincreasing_multi}
Let $\boldsymbol{M} \in \mathbb{R}^d$ be a random vector with $M_i>0\; \forall i=1,2,\ldots,d$ and $\epsilon > 0$. For every continuous function $f: \mathbb{R}^d \rightarrow \mathbb{R}$ such that $f(\boldsymbol{m})$ is increasing and $f(\boldsymbol{m}) \geq 0$, for $m_i > 0\; \forall i \in \left\lbrace 1,2,\ldots, d \right\rbrace$, where $\boldsymbol{m} = \left(m_1,m_2,\ldots,m_d\right)$,
\begin{equation}
\nonumber \EE[f(\boldsymbol{M}) | M_i < \epsilon\; \forall i=1,2,\ldots,d] \leq \EE[f(\boldsymbol{M})].
\end{equation}
\end{lemma}
\vspace{0.1in}
{\raggedright \textit{Proof}{\rm of \textbf{Theorem \ref{Theorem_i.n.i.d}}.}}
It has already been shown in the outline of the proof that the triangle inequality yields
\begin{align}
\nonumber & \left|\EE\left[h\left(\sqrt{n}\left[\bar{I}_n(\boldsymbol{\theta_0})\right]^{\frac{1}{2}}\left(\boldsymbol{\hat{\theta}_n(X)} - \boldsymbol{\theta_0}\right)\right)\right]- \EE[h(\boldsymbol{Z})]\right| \leq \eqref{first_term_regression} + \eqref{second_term_regression}.
\end{align}
{\textbf{Step 1: Upper bound for} \eqref{first_term_regression}}. First, $\nabla (\ell(\boldsymbol{\theta_0};\boldsymbol{x})) = \sum_{i=1}^{n}\nabla \left(\log(f_i(\boldsymbol{x_i}|\boldsymbol{\theta_0}))\right)$ due to independence. With $\tilde{V}$ as in \eqref{cm}, the results of Theorem 2.1 of \cite{ReinertRollin} will be used for
\begin{equation}
\label{W_i.n.i.d.}
\boldsymbol{W} = \frac{1}{\sqrt{n}}\tilde{V}\sum_{i=1}^{n}\nabla(\log(f_i(\boldsymbol{X_i}|\boldsymbol{\theta_0}))) = \left(W_1, W_2, \ldots, W_d\right)^{\intercal} \in \mathbb{R}^{d \times 1}.
\end{equation}
From \eqref{W_i.n.i.d.} we have that for all $k \in \left\lbrace 1,2, \ldots, d\right\rbrace$, $W_k = \sum_{i=1}^{n} \xi_{ik}$, with $\xi_{ik}$ as in \eqref{cm}. From the regularity conditions, $\EE\left[\nabla(\ell(\boldsymbol{\theta_0};\boldsymbol{X}))\right] = \boldsymbol{0}$ and thus $\EE[\boldsymbol{W}] = \boldsymbol{0}$. Also, $\bar{I}_n(\boldsymbol{\theta_0})$ is symmetric. Therefore, $\tilde{V}$ is also symmetric. Using the regularity conditions we know that\linebreak $\sum_{i=1}^{n}{\rm Cov}\left[\nabla\left(\log\left(f_i(\boldsymbol{X_i}|\boldsymbol{\theta_0})\right)\right)\right] = n\bar{I}_n(\boldsymbol{\theta_0})$ and basic calculations show that ${\rm Cov}[\boldsymbol{W}] = I_{d \times d}$. Since $\EE[\boldsymbol{W}] = \boldsymbol{0}$ and $\EE\left[\boldsymbol{W}\boldsymbol{W}^{\intercal}\right] = I_{d\times d}$, the first assumption of Theorem 2.1 from \cite{ReinertRollin} is satisfied. This theorem also assumes that $\exists \boldsymbol{W'}$ such that $(\boldsymbol{W}, \boldsymbol{W'})$ is an exchangeable pair meaning that $(\boldsymbol{W},\boldsymbol{W'})\stackrel{\mathclap{\normalfont\mbox{d}}}{=\joinrel=}(\boldsymbol{W'},\boldsymbol{W})$, where $\stackrel{\mathclap{\normalfont\mbox{d}}}{=\joinrel=}$ denotes equality in distribution. In addition, it is assumed that
\begin{equation}
\label{lSteinpair}
\EE\left[\boldsymbol{W'} - \boldsymbol{W}|\boldsymbol{W}\right] = -\Lambda\boldsymbol{W} + \boldsymbol{R}
\end{equation}
for an invertible $d \times d$ matrix $\Lambda$ and a $\sigma(\boldsymbol{W})$-measurable random vector $\boldsymbol{R}$. To define $\boldsymbol{W'}$ in our case such that \eqref{lSteinpair} is satisfied, let $\left\lbrace \boldsymbol{X_i'}, i=1,2,\ldots,n \right\rbrace$ be an independent copy of $\left\lbrace \boldsymbol{X_i}, i=1,2,\ldots,n\right\rbrace$ and let the index $I \in \left\lbrace 1,2,\ldots,n\right\rbrace$ follow the uniform distribution on $\left\lbrace 1,2,\ldots,n\right\rbrace$, independently of the set $\left\lbrace \boldsymbol{X_i},\boldsymbol{X_i'}, i=1,2,\ldots,n \right\rbrace$. Let $$\xi_{ik}' = \frac{1}{\sqrt{n}}\sum_{j=1}^{d}\tilde{V}_{kj}\frac{\partial}{\partial \theta_j}\log(f_i(\boldsymbol{X_i'}|\boldsymbol{\theta_0}))$$ and $$W_k'=W_k-\xi_{Ik}+\xi_{Ik}',\; \forall k \in \left\lbrace 1,2,\ldots,d\right\rbrace,$$
with $\EE\left[W_{k}' - W_{k}|\boldsymbol{W}\right] = \EE\left[\xi_{Ik}' - \xi_{Ik}|\boldsymbol{W}\right] = -\EE\left[\xi_{Ik}|\boldsymbol{W}\right] = -\frac{1}{n}\sum_{i=1}^{n}\EE\left[\xi_{ik}|\boldsymbol{W}\right] = -\frac{W_{k}}{n}.$ Hence \eqref{lSteinpair} is satisfied with $\Lambda = \frac{1}{n}I_{d \times d}$ and $\boldsymbol{R} = \boldsymbol{0}$. Therefore, Theorem 2.1 from \cite{ReinertRollin} gives in our case that
\begin{align}
\label{A1}
\left|\EE[h(\boldsymbol{W})] - \EE[h(\boldsymbol{Z})]\right| \leq & n\left(\frac{\|h\|_2}{4}\sum_{i=1}^{d}\sum_{j=1}^{d}\left[{\rm Var}\left[\EE\left[\left(W_i'-W_i\right)\left(W_j'-W_j\right)|\boldsymbol{W}\right]\right]\right]^{\frac{1}{2}}\right)\\
\label{A2}
& + n\left(\frac{\|h\|_3}{12}\sum_{i=1}^{d}\sum_{j=1}^{d}\sum_{k=1}^{d}\EE\left|\left(W_{i}'-W_i\right)\left(W_{j}' - W_j\right)\left(W_{k}'-W_{k}\right)\right|\right).
\end{align}
To bound the variance of the conditional expectations in \eqref{A1}, let $\mathscr{A} = \sigma\left(\boldsymbol{X_1}, \boldsymbol{X_2}, \ldots, \boldsymbol{X_n}\right)$. Since $\sigma(\boldsymbol{W}) \subset \mathscr{A}$, for any random variable $Y$, we have that ${\rm Var}\left[\EE[Y|\boldsymbol{W}]\right] \leq {\rm Var}\left[\EE\left[Y|\mathscr{A}\right]\right]$. Then,
\begin{align}
\label{Amidstep1}
\nonumber & \eqref{A1} \leq n\frac{\|h\|_2}{4}\left\lbrace\sum_{j=1}^{d}\sqrt{{\rm Var}\left[\EE\left[(\xi_{Ij}'-\xi_{Ij})^2|\mathscr{A}\right]\right]}\right.\\
& \left.\qquad\qquad\qquad + 2\sum_{k=1}^{d-1}\sum_{j=k+1}^{d}\sqrt{{\rm Var}\left[\EE\left[\left(\xi_{Ik}'-\xi_{Ik}\right)\left(\xi_{Ij}'-\xi_{Ij}\right)|\mathscr{A}\right]\right]}\right\rbrace.
\end{align}
Since $\left\lbrace \boldsymbol{X_i'}, i=1,2,\ldots,n \right\rbrace$ is an independent copy of $\left\lbrace \boldsymbol{X_i}, i=1,2,\ldots,n\right\rbrace$ and $\xi_{ik}'$ is independent of $\mathscr{A}$,
\begin{align}
\label{Amidstep2}
\nonumber & \eqref{Amidstep1} = n\frac{\|h\|_2}{4}\left\lbrace\vphantom{(\left(\sup_{\theta:|\theta-\theta_0|\leq\epsilon}\left|l^{(3)}(\theta;\boldsymbol{X})\right|\right)^2}\sum_{j=1}^{d}\left[{\rm Var}\left[\EE\left[(\xi_{Ij}')^2\right] -2\EE[\xi_{Ij}']\EE[\xi_{Ij}|\mathscr{A}] + \EE\left[\xi_{Ij}^2|\mathscr{A}\right]\right]\right]^{\frac{1}{2}}\right.\\
&\left. + 2\sum_{k=1}^{d-1}\sum_{j=k+1}^{d}\left[{\rm Var}\left[\EE\left[\xi_{Ik}'\xi_{Ij}'\right] - \EE\left[\xi_{Ij}'\right]\EE\left[\xi_{Ik}|\mathscr{A}\right] - \EE\left[\xi_{Ik}'\right]\EE\left[\xi_{Ij}|\mathscr{A}\right] + \EE\left[\xi_{Ik}\xi_{Ij}|\mathscr{A}\right]\right]\right]^{\frac{1}{2}}\vphantom{(\left(\sup_{\theta:|\theta-\theta_0|\leq\epsilon}\left|l^{(3)}(\theta;\boldsymbol{X})\right|\right)^2}\right\rbrace.
\end{align}
Using that $\EE\left[\xi_{ik}'\right] = 0$,
\begin{align}
%\nonumber & \eqref{Amidstep2} = n\frac{\|h\|_2}{4}\left\lbrace\vphantom{(\left(\sup_{\theta:|\theta-\theta_0|\leq\epsilon}\left|l^{(3)}(\theta;\boldsymbol{X})\right|\right)^2}\sum_{j=1}^{d}\left[{\rm Var}\left(\EE\left[(\xi_{Ij}')^2\right] + \EE\left[\xi_{Ij}^2|\mathscr{A}\right]\right)\right]^{\frac{1}{2}}\right.\\
%\nonumber & \quad\qquad\qquad\qquad\left. + 2\sum_{k=1}^{d-1}\sum_{j=k+1}^{d}\left[{\rm Var}\left(\EE\left[\xi_{Ik}'\xi_{Ij}'\right] + \EE\left[\xi_{Ik}\xi_{Ij}|\mathscr{A}\right]\right)\right]^{\frac{1}{2}}\vphantom{(\left(\sup_{\theta:|\theta-\theta_0|\leq\epsilon}\left|l^{(3)}(\theta;\boldsymbol{X})\right|\right)^2}\right\rbrace\\
\nonumber & \eqref{Amidstep2} = n\frac{\|h\|_2}{4}\left\lbrace\vphantom{(\left(\sup_{\theta:|\theta-\theta_0|\leq\epsilon}\left|l^{(3)}(\theta;\boldsymbol{X})\right|\right)^2}\sum_{j=1}^{d}\left[\frac{1}{n^2}{\rm Var}\left[\sum_{i=1}^{n}\EE\left[\xi_{ij}^2|\mathscr{A}\right]\right]\right]^{\frac{1}{2}} + 2\sum_{k=1}^{d-1}\sum_{j=k+1}^{d}\left[\frac{1}{n^2}{\rm Var}\left[\sum_{i=1}^{n}\EE\left[\xi_{ik}\xi_{ij}|\mathscr{A}\right]\right]\right]^{\frac{1}{2}}\vphantom{(\left(\sup_{\theta:|\theta-\theta_0|\leq\epsilon}\left|l^{(3)}(\theta;\boldsymbol{X})\right|\right)^2}\right\rbrace\\
\nonumber & = \frac{\|h\|_2}{4}\left\lbrace\sum_{j=1}^{d}\left[{\rm Var}\left[\sum_{i=1}^{n}\xi_{ij}^2\right]\right]^{\frac{1}{2}} + 2\sum_{k=1}^{d-1}\sum_{j=k+1}^{d}\left[{\rm Var}\left[\sum_{i=1}^{n}\xi_{ik}\xi_{ij}\right]\right]^{\frac{1}{2}}\right\rbrace\\
\nonumber & = \frac{\|h\|_2}{\sqrt{n}}K_2(\boldsymbol{\theta_0}),
\end{align}
with $K_2(\boldsymbol{\theta_0})$ defined in \eqref{remainder_term_R2}. For \eqref{A2}, using again the definition of $\xi_{ik}$ in \eqref{cm}, after basic calculations we obtain that
\begin{equation}
\nonumber \eqref{A2} \leq \frac{\|h\|_3}{\sqrt{n}}K_3(\boldsymbol{\theta_0}),
\end{equation}
with $K_3(\boldsymbol{\theta_0})$ as in \eqref{remainder_term_R3}. Therefore,
\begin{equation}
\label{bound_first_regression}
\eqref{first_term_regression} \leq \frac{\|h\|_2}{\sqrt{n}}K_2(\boldsymbol{\theta_0}) +\frac{\|h\|_3}{\sqrt{n}}K_3(\boldsymbol{\theta_0}).
\end{equation}
{\textbf{Step 2: Upper bound for} \eqref{second_term_regression}}. With $\tilde{V}$ as in \eqref{cm}, for ease of presentation let us denote by
\begin{align}
\label{notationmultiT1T2}
\nonumber & \boldsymbol{R_1(\theta_0;x)} = \frac{1}{2\sqrt{n}}\tilde{V}\sum_{j=1}^{d}\sum_{q=1}^{d}Q_jQ_q\left(\nabla\left(\frac{\partial^2}{\partial\theta_j\partial\theta_q}\ell(\boldsymbol{\theta};\boldsymbol{x})\Big|_{\substack{\boldsymbol{\theta} = \boldsymbol{\theta_0^{*}}}}\right)\right)\\
\nonumber & T_1 = T_1(\boldsymbol{\theta_0};\boldsymbol{X},h) := h\left(\sqrt{n}\left[\bar{I}_n(\boldsymbol{\theta_0})\right]^{\frac{1}{2}}\left(\boldsymbol{\hat{\theta}_n(X)} - \boldsymbol{\theta_0}\right)\right) - h\left(\vphantom{(\left(\sup_{\theta:|\theta-\theta_0|\leq\epsilon}\left|l^{(3)}(\theta;\boldsymbol{X})\right|\right)^2}\frac{1}{\sqrt{n}}\tilde{V}\left(\nabla (\ell(\boldsymbol{\theta_0};\boldsymbol{x}))\right) + \boldsymbol{R_1(\theta_0;X)}\right)\\
& T_2 = T_2(\boldsymbol{\theta_0};\boldsymbol{X},h) := h\left(\vphantom{(\left(\sup_{\theta:|\theta-\theta_0|\leq\epsilon}\left|l^{(3)}(\theta;\boldsymbol{X})\right|\right)^2}\frac{1}{\sqrt{n}}\tilde{V}\left(\nabla (\ell(\boldsymbol{\theta_0};\boldsymbol{x}))\right) + \boldsymbol{R_1(\theta_0;x)}\right) -  h\left(\frac{1}{\sqrt{n}}\tilde{V}\left(\nabla\left(\ell(\boldsymbol{\theta_0};\boldsymbol{X})\right)\right)\right).
\end{align}
Using the above notation and the triangle inequality
\begin{equation}
\nonumber \eqref{second_term_regression} = \left|\EE\left[T_1 + T_2\right]\right| \leq \EE|T_1| + \EE|T_2|.
\end{equation}
With $A_{[j]}$ the $j^{{\rm th}}$ row of a matrix $A$, a first order multivariate Taylor expansion gives that
\begin{align}
\nonumber & \left|T_1\right| \leq \|h\|_1\left|\vphantom{(\left(\sup_{\theta:|\theta-\theta_0|\leq\epsilon}\left|l^{(3)}(\theta;\boldsymbol{X})\right|\right)^2}\sum_{j=1}^{d}\left(\vphantom{(\left(\sup_{\theta:|\theta-\theta_0|\leq\epsilon}\left|l^{(3)}(\theta;\boldsymbol{X})\right|\right)^2}\sqrt{n}\left[\left[\bar{I}_n(\boldsymbol{\theta_0})\right]^{\frac{1}{2}}\right]_{[j]}(\boldsymbol{\hat{\theta}_n(X)} - \boldsymbol{\theta_0})-\frac{1}{\sqrt{n}}\tilde{V}_{[j]}\nabla\left(\ell(\boldsymbol{\theta_0};\boldsymbol{X})\right)\right.\right.\\
\nonumber & \qquad\quad\qquad\quad\left.\left.- \frac{1}{2\sqrt{n}}\tilde{V}_{[j]}\left\lbrace\vphantom{(\left(\sup_{\theta:|\theta-\theta_0|\leq\epsilon}\left|l^{(3)}(\theta;\boldsymbol{X})\right|\right)^2}\sum_{k=1}^{d}\sum_{q=1}^{d}Q_kQ_q\left(\nabla\left(\frac{\partial^2}{\partial\theta_k\partial\theta_q}\ell(\boldsymbol{\theta};\boldsymbol{x})\Big|_{\substack{\boldsymbol{\theta} = \boldsymbol{\theta_0^{*}}}}\right)\right)\vphantom{(\left(\sup_{\theta:|\theta-\theta_0|\leq\epsilon}\left|l^{(3)}(\theta;\boldsymbol{X})\right|\right)^2}\right\rbrace\vphantom{(\left(\sup_{\theta:|\theta-\theta_0|\leq\epsilon}\left|l^{(3)}(\theta;\boldsymbol{X})\right|\right)^2}\right)\vphantom{(\left(\sup_{\theta:|\theta-\theta_0|\leq\epsilon}\left|l^{(3)}(\theta;\boldsymbol{X})\right|\right)^2}\right|.
\end{align}
Using \eqref{multiresultTaylor} component-wise and the Cauchy-Schwarz inequality, we have that
\begin{align}
\label{boundforT1first}
& \EE[T_1] \leq\frac{\|h\|_1}{\sqrt{n}} \sum_{k=1}^{d}\sum_{l=1}^{d}\left|\tilde{V}_{lk}\right|\sum_{j=1}^{d}\sqrt{\EE\left[Q_j^2\right]\EE\left[\left(\frac{\partial^2}{\partial \theta_j\partial\theta_k}\ell(\boldsymbol{\theta_0};\boldsymbol{X}) + n[\bar{I}_n(\boldsymbol{\theta_0})]_{kj}\right)^2\right]}.
%\frac{\|h\|_1}{\sqrt{n}} \sum_{k=1}^{d}\sum_{l=1}^{d}\left|\tilde{V}_{lk}\right|\sum_{j=1}^{d}\EE\left|\vphantom{(\left(\sup_{\theta:|\theta-\theta_0|\leq\epsilon}\left|l^{(3)}(\theta;\boldsymbol{X})\right|\right)^2}Q_j\left(\frac{\partial^2}{\partial \theta_j\partial\theta_k}\ell(\boldsymbol{\theta_0};\boldsymbol{X}) + n[\bar{I}_n(\boldsymbol{\theta_0})]_{kj}\right)\vphantom{(\left(\sup_{\theta:|\theta-\theta_0|\leq\epsilon}\left|l^{(3)}(\theta;\boldsymbol{X})\right|\right)^2}\right|
\end{align}
To bound now $\EE\left|T_2\right|$, with $T_2$ as in \eqref{notationmultiT1T2}, we need to take into account that $\frac{\partial^3}{\partial\theta_k\partial\theta_q\partial\theta_j}\ell(\boldsymbol{\theta};\boldsymbol{x})\Big|_{\substack{\boldsymbol{\theta} = \boldsymbol{\theta_0^{*}}}}$ is in general not uniformly bounded. For $\epsilon>0$, the law of total expectation and Markov's inequality yield
\begin{align}
\label{chebyshevmultiT2}
\EE\left|T_2\right| &\leq 2\|h\|\Prob\left(\left|Q_{(m)}\right|\geq\epsilon\right) + \EE\left[|T_2|\middle|\left|Q_{(m)}\right| < \epsilon\right] \leq \frac{2\|h\|}{\epsilon^2}\EE\left[\sum_{j=1}^{d}Q_j^2\right] + \EE\left[|T_2|\middle|\left|Q_{(m)}\right| < \epsilon\right],
\end{align}
with $Q_{(m)}$ as in \eqref{cm}. To bound $\EE\left[|T_2|\middle|\left|Q_{(m)}\right| < \epsilon\right]$, a first-order Taylor expansion  and \eqref{multiresultTaylor} yield
\begin{align}
\label{T2multijustabound}
\left|T_2\right| &\leq \frac{\|h\|_1}{2\sqrt{n}}\sum_{k=1}^{d}\sum_{l=1}^{d}\left|\tilde{V}_{lk}\right|\sum_{j=1}^{d}\sum_{v=1}^{d}\left|Q_jQ_v\frac{\partial^3}{\partial\theta_{k}\partial\theta_{j}\partial\theta_{v}}\ell(\boldsymbol{\theta};\boldsymbol{X})\Big|_{\substack{\boldsymbol{\theta} = \boldsymbol{\theta_0^{*}}}}\right|.
\end{align}
Therefore, from \eqref{chebyshevmultiT2} and \eqref{T2multijustabound} we have that
\begin{align}
\nonumber & \EE|T_2| \leq \frac{2\|h\|}{\epsilon^2}\EE\left[\sum_{j=1}^{d}Q_j^2\right] +\frac{\|h\|_1}{2\sqrt{n}}\sum_{k=1}^{d}\sum_{l=1}^{d}\left|\tilde{V}_{lk}\right|\EE\left[\vphantom{(\left(\sup_{\theta:|\theta-\theta_0|\leq\epsilon}\left|l^{(3)}(\theta;\boldsymbol{X})\right|\right)^2}\sum_{j=1}^{d}\sum_{v=1}^{d}\left|\vphantom{(\left(\sup_{\theta:|\theta-\theta_0|\leq\epsilon}\left|l^{(3)}(\theta;\boldsymbol{X})\right|\right)^2}Q_jQ_v\frac{\partial^3}{\partial\theta_{k}\partial\theta_{j}\partial\theta_{v}}\ell(\boldsymbol{\theta};\boldsymbol{X})\Big|_{\substack{\boldsymbol{\theta} = \boldsymbol{\theta_0^{*}}}}\vphantom{(\left(\sup_{\theta:|\theta-\theta_0|\leq\epsilon}\left|l^{(3)}(\theta;\boldsymbol{X})\right|\right)^2}\right|\middle|\left|Q_{(m)}\right| < \epsilon\vphantom{(\left(\sup_{\theta:|\theta-\theta_0|\leq\epsilon}\left|l^{(3)}(\theta;\boldsymbol{X})\right|\right)^2}\right].
\end{align}
The Cauchy-Schwarz inequality and Lemma \ref{Lemmaincreasing_multi} yield
\begin{align}
\label{T2expfin}
\nonumber & \EE|T_2|  \leq \frac{2\|h\|}{\epsilon^2}\EE\left[\sum_{j=1}^{d}Q_j^2\right]\\
& \; +  \frac{\|h\|_1}{2\sqrt{n}}\left\lbrace\vphantom{\left[\EE\left(\left(\sup_{\boldsymbol{\theta}:\left|\theta_m - \theta_{0,m}\right|\leq \epsilon}\left|\frac{\partial^3}{\partial \theta_{k}\partial \theta_{j}\partial \theta_{i}}\ell(\boldsymbol{\theta};\boldsymbol{X})\right|\right)^2\middle|\left|Q_{(m)}\right|\leq \epsilon \right)\right]^{\frac{1}{2}}}\sum_{k=1}^{d}\sum_{l=1}^{d}\left|\tilde{V}_{lk}\right|\sum_{j=1}^{d}\sum_{v=1}^{d}\left[\EE\left[Q_j^2Q_v^2\right]\right]^{\frac{1}{2}}\left[\EE\left[\left(M_{kjv}(\boldsymbol{X})\right)^2\middle|\left|Q_{(m)}\right| < \epsilon \right]\right]^{\frac{1}{2}}\vphantom{\left[\EE\left(\left(\sup_{\boldsymbol{\theta}:\left|\theta_m - \theta_{0,m}\right|\leq \epsilon}\left|\frac{\partial^3}{\partial \theta_{k}\partial \theta_{j}\partial \theta_{i}}\ell(\boldsymbol{\theta};\boldsymbol{X})\right|\right)^2\middle|\left|Q_{(m)}\right|\leq \epsilon \right)\right]^{\frac{1}{2}}}\right\rbrace.
%\frac{2\|h\|}{\epsilon^2}\EE\left(\sum_{j=1}^{d}\left(\hat{\theta}_n(\boldsymbol{X})_j - \theta_{0,j}\right)^2\right)\\
%\nonumber & +  \frac{\|h\|_1}{2\sqrt{n}}\left\lbrace\vphantom{\left[\EE\left(\left(\sup_{\boldsymbol{\theta}:\left|\theta_m - \theta_{0,m}\right|\leq \epsilon}\left|\frac{\partial^3}{\partial \theta_{k}\partial \theta_{j}\partial \theta_{i}}\ell(\boldsymbol{\theta};\boldsymbol{X})\right|\right)^2\middle|\left|Q_{(m)}\right|\leq \epsilon \right)\right]^{\frac{1}{2}}}\sum_{k=1}^{d}\sum_{l=1}^{d}\left|\left[\left[\bar{I}_n(\boldsymbol{\theta_0})\right]^{-\frac{1}{2}}\right]_{lk}\right|\sum_{j=1}^{d}\sum_{v=1}^{d}\left[\EE\left(\left(\hat{\theta}_n(\boldsymbol{X})_j - \theta_{0,j}\right)^2\left(\hat{\theta}_n(\boldsymbol{X})_v - \theta_{0,v}\right)^2\right)\right]^{\frac{1}{2}}\right.\\
%\nonumber &\qquad\quad \times \left.\left[\EE\left(\left(\sup_{\substack{\boldsymbol{\theta}:\left|\theta_m - \theta_{0,m}\right|< \epsilon\\ \forall m \in \left\lbrace 1,2,\ldots, d\right\rbrace}}\left|\frac{\partial^3}{\partial\theta_k\partial\theta_j\partial\theta_v}\ell(\boldsymbol{\theta};\boldsymbol{X})\right|\right)^2\middle|\left|Q_{(m)}\right| < \epsilon \right)\right]^{\frac{1}{2}}\vphantom{\left[\EE\left(\left(\sup_{\boldsymbol{\theta}:\left|\theta_m - \theta_{0,m}\right|\leq \epsilon}\left|\frac{\partial^3}{\partial \theta_{k}\partial \theta_{j}\partial \theta_{i}}\ell(\boldsymbol{\theta};\boldsymbol{X})\right|\right)^2\middle|\left|Q_{(m)}\right|\leq \epsilon \right)\right]^{\frac{1}{2}}}\right\rbrace\\
\end{align}
Therefore, from \eqref{boundforT1first} and \eqref{T2expfin} we obtain that
\begin{equation}
\label{boundmultivariateB}
\eqref{second_term_regression} \leq \frac{2\|h\|}{\epsilon^2}\EE\left[\sum_{j=1}^{d}Q_j^2\right] + \frac{\|h\|_1}{\sqrt{n}}K_1(\boldsymbol{\theta_0}),
\end{equation}
where $K_1(\boldsymbol{\theta_0})$ is as in \eqref{remainder_term_R1}. Using now \eqref{bound_first_regression} and \eqref{boundmultivariateB} we obtain the assertion.$\qquad\qquad\qquad\quad\blacksquare$

\section*{Appendix: Proofs of Lemma \ref{Lemmaincreasing_multi} and of Corollaries \ref{Corollary_multi_normal_bound} and \ref{Cor_multi_Beta}}

{\raggedright {\textit{Proof}}{\rm\;of\;\textbf{Lemma \ref{Lemmaincreasing_multi}}.}}  Let $k \in \left\lbrace 1,2,\ldots,d \right\rbrace$. We set $M_{d+1}=0$. It will be shown that for $k=1,2,\ldots,d$ we have that
\begin{equation}
\nonumber \EE[f(\boldsymbol{M}) | M_i < \epsilon\;,\; i=k,\ldots,d] \leq \EE[f(\boldsymbol{M}) | M_i < \epsilon\;,\; i=k+1,\ldots,d].
\end{equation}
From the law of total expectation,
\begin{align}
\nonumber & \EE[f(\boldsymbol{M}) | M_i < \epsilon\;,\; i=k+1,\ldots,d]\\
\nonumber & = \EE[f(\boldsymbol{M}) | M_i < \epsilon\;,\; i=k,\ldots,d]\Prob\left[M_{k} < \epsilon|M_{i} < \epsilon\;,\;i=k+1,\ldots,d\right]\\
\nonumber & \; + \EE[f(\boldsymbol{M}) | M_i < \epsilon\;,\; i=k+1,\ldots,d,\; M_{k}\geq\epsilon]\Prob\left[M_{k}\geq\epsilon|M_{i} < \epsilon\;,\;i=k+1,\ldots,d\right].
\end{align}
Using that
\begin{equation}
\nonumber \Prob\left[M_{k} < \epsilon|M_{i} < \epsilon\;,\;i=k+1,\ldots,d\right] = 1 - \Prob\left[M_{k} \geq \epsilon|M_{i} < \epsilon\;,\;i=k+1,\ldots,d\right]
\end{equation}
yields
\begin{align}
\label{midLemmaincreasingproof}
\nonumber & \EE[f(\boldsymbol{M}) | M_i < \epsilon\;,\; i=k+1,\ldots,d] = \EE[f(\boldsymbol{M}) | M_i < \epsilon\;,\; i=k,\ldots,d]\\
\nonumber & +\Prob\left[M_{k} \geq \epsilon|M_{i} < \epsilon\;,\;i=k+1,\ldots,d\right]\left\lbrace\vphantom{(\left(\sup_{\theta:|\theta-\theta_0|\leq\epsilon}\left|l^{(3)}(\theta;\boldsymbol{X})\right|\right)^2}\EE[f(\boldsymbol{M}) | M_i < \epsilon\;,\; i=k+1,\ldots,d,\; M_{k}\geq\epsilon]\right.\\
& \qquad\qquad\qquad\qquad\qquad\qquad\qquad\qquad\quad\left.- \EE[f(\boldsymbol{M}) | M_i < \epsilon\;,\; i=k,\ldots,d]\vphantom{(\left(\sup_{\theta:|\theta-\theta_0|\leq\epsilon}\left|l^{(3)}(\theta;\boldsymbol{X})\right|\right)^2}\right\rbrace.
\end{align}
Since $f(\boldsymbol{m})$ is an increasing function,
\begin{align}
\nonumber \EE[f(\boldsymbol{M}) | M_i < \epsilon\;,\; i=k+1,\ldots,d,\; M_{k}\geq\epsilon] - \EE[f(\boldsymbol{M}) | M_i < \epsilon\;,\; i=k,\ldots,d]\geq 0.
\end{align}
Applying this to \eqref{midLemmaincreasingproof} gives that
\begin{equation}
\nonumber \EE[f(\boldsymbol{M}) | M_i < \epsilon\;,\; i=k,\ldots,d] \leq \EE[f(\boldsymbol{M}) | M_i < \epsilon\;,\; i=k+1,\ldots,d].
\end{equation}
A simple iteration over $k$ gives that
\begin{equation}
\nonumber \EE[f(\boldsymbol{M}) | M_i < \epsilon\; \forall i=1,2,\ldots,d] \leq \EE[f(\boldsymbol{M})],
\end{equation}
which is the result of the lemma. $\quad\quad\qquad\qquad\qquad\qquad\qquad\qquad\qquad\qquad\qquad\qquad\qquad\quad\;\;\blacksquare$
\vspace{0.1in}
\\
{\raggedright \textit{Proof}{\rm of \textbf{Corollary \ref{Corollary_multi_normal_bound}}.}}
For one random variable, the first and second-order partial derivatives of the logarithm of the normal density function are
\begin{align}
\label{derivatives_multi_normal}
\nonumber & \frac{\partial}{\partial \eta_1}\log f(x_1|\boldsymbol{\eta_0}) = -x_1^2 + \frac{1}{2\eta_1} + \frac{\eta_2^2}{4\eta_1^2},\qquad\qquad \frac{\partial}{\partial \eta_2}\log f(x_1|\boldsymbol{\eta_0}) = x_1 - \frac{\eta_2}{2\eta_1},\\
\nonumber &\frac{\partial^2}{\partial\eta_1^2}\log f(x_1|\boldsymbol{\eta_0}) = -\left(\frac{1}{2\eta_1^2} + \frac{\eta_2^2}{2\eta_1^3}\right), \qquad\qquad\; \frac{\partial^2}{\partial\eta_2^2} \log f(x_1|\boldsymbol{\eta_0}) = -\frac{1}{2\eta_1},\\
&\frac{\partial^2}{\partial\eta_1\partial\eta_2} \log f(x_1|\boldsymbol{\eta_0}) = \frac{\partial^2}{\partial\eta_2\partial\eta_1} \log f(x_1|\boldsymbol{\eta_0}) = \frac{\eta_2}{2\eta_1^2}.
\end{align}
Hence, the expected Fisher Information matrix for one random variable is
\begin{equation}
\label{multi_normal_FISHER1}
I(\boldsymbol{\theta_0}) = \frac{1}{2\eta_1}\begin{pmatrix}
\frac{1}{\eta_1} + \frac{\eta_2^2}{\eta_1^2} & -\frac{\eta_2}{\eta_1}\\
-\frac{\eta_2}{\eta_1} & 1
\end{pmatrix},
\end{equation}
and after simple calculations we obtain that
\begin{equation}
\label{multi_normal_FISHER}
\nonumber \left[I(\boldsymbol{\theta_0})\right]^{-\frac{1}{2}} = \tilde{V} = \sqrt{\frac{2}{\alpha}}\begin{pmatrix}
\eta_1^{\frac{3}{2}}\left(1+\sqrt{\eta_1}\right) & \eta_1\eta_2\\
\eta_1\eta_2 & \eta_1\left(1+\sqrt{\eta_1}\right)+\eta_2^2
\end{pmatrix},
\end{equation}
where $\alpha = \eta_1\left(1+\sqrt{\eta_1}\right)^2 + \eta_2^2$ as defined in Corollary \ref{Corollary_multi_normal_bound}. We bound the terms in Theorem \ref{Theorem_i.n.i.d} in order of appearance. The term $K_1(\boldsymbol{\eta_0})$ is given in \eqref{remainder_term_R1} and the first quantity of $K_1(\boldsymbol{\eta_0})$ vanishes due to the fact that
\begin{equation}
\label{Tkjzero}
\EE\left[T_{kj}^2\right] = 0, \quad \forall k,j \in \left\lbrace 1,2\right\rbrace.
\end{equation}
This comes from the definition of $T_{kj}$ in \eqref{cm} and the results of \eqref{derivatives_multi_normal} and \eqref{multi_normal_FISHER1}. For the second term of $K_1(\boldsymbol{\eta_0})$, we note that ${\rm Cov}\left[\bar{X}, \frac{1}{n}\sum_{i=1}^{n}\left(X_i - \bar{X}\right)^2\right] = 0$ \citep{Casella}[p.218] and simple calculations lead to
\begin{align}
\label{mixed_MSE}
\nonumber & \EE\left[Q_1^2Q_2^2\right] < \frac{2n\eta_1^3(2n+63) + 3\eta_1^2\eta_2^2\left(4n^2+172n+315\right)}{(n-5)^2(n-9)^2},\\
\nonumber & \EE\left[Q_1^4\right] < \frac{\eta_1^4\left(12n^2+516n+945\right)}{(n-5)^2(n-9)^2}\\
& \EE\left[Q_2^4\right] < \frac{12n^2\left(\eta_1+\eta_2^2\right)^2+12n\eta_2^2\left(43\eta_2^2+63\eta_1\right)+945\eta_2^4}{(n-5)^2(n-9)^2},
\end{align}
where $Q_1$ and $Q_2$ are defined in \eqref{cm}. In addition, for $M_{kjl}(\boldsymbol{x})$ and $0 < \epsilon = \epsilon(\boldsymbol{\eta_0})$ as in the condition (R.C.3), simple calculations and \eqref{derivatives_multi_normal} yield for $m =1,2$
\begin{align}
\label{Mjklterms}
\nonumber & \underset{\boldsymbol{\theta}:\left|\theta_m - \eta_{m}\right| < \epsilon}{\sup}\left|\frac{\partial^3}{\partial \theta_1^3}\ell(\boldsymbol{\theta};\boldsymbol{X})\right| = \underset{\boldsymbol{\theta}:\left|\theta_m - \eta_{m}\right| < \epsilon}{\sup}\left|\frac{n}{\theta_1^3} + \frac{3n\theta_2^2}{2\theta_1^4}\right|< \frac{n}{\left(\eta_1-\epsilon\right)^3}\left(1+\frac{3(\eta_2+\epsilon)^2}{2(\eta_1-\epsilon)}\right)=: M_{111}(\boldsymbol{x}),\\
\nonumber & \underset{\boldsymbol{\theta}:\left|\theta_m - \eta_{m}\right| < \epsilon}{\sup}\left|\frac{\partial^3}{\partial \theta_2^3}\ell(\boldsymbol{\theta};\boldsymbol{X})\right| = 0 =: M_{222}(\boldsymbol{x}),\\
\nonumber & \underset{\boldsymbol{\theta}:\left|\theta_m - \eta_{m}\right| < \epsilon}{\sup}\left|\frac{\partial^3}{\partial \theta_1^2\theta_2}\ell(\boldsymbol{\theta};\boldsymbol{X})\right| = \left|-\frac{n\eta_2}{\eta_1^3}\right| <  \frac{n(\eta_2+\epsilon)}{(\eta_1-\epsilon)^3} =: M_{112}(\boldsymbol{x}),\\
& \underset{\boldsymbol{\theta}:\left|\theta_m - \eta_{m}\right| < \epsilon}{\sup}\left|\frac{\partial^3}{\partial \theta_1\theta_2^2}\ell(\boldsymbol{\theta};\boldsymbol{X})\right| = \left|\frac{n}{2\eta_1^2}\right| < \frac{n}{2(\eta_1-\epsilon)^2} =: M_{221}(\boldsymbol{x}). 
\end{align}
For the choice of $\epsilon = \epsilon_0$ as in (R.C.3), \eqref{Mjklterms} requires that $0< \epsilon<\eta_1$. There is a trade-off on its choice for the fourth term of the bound in \eqref{final_bound_regression} and the results in \eqref{Mjklterms}. This is because the last term of the general bound in \eqref{final_bound_regression} is divided by $\epsilon^2$ indicating that we should choose $\epsilon$ away from zero. However the terms in \eqref{Mjklterms} have powers of $\eta_1 - \epsilon$ on the denominator and it would be reasonable for $\epsilon$ to be close to zero and away from $\eta_1$. An optimisation process with respect to $\epsilon$ becomes quite tedious and therefore we choose $\epsilon$ to be the midpoint of $(0,\eta_1)$, which is sufficiently away from both zero and $\eta_1$ and also behaves very well. Using this value of $\epsilon$ and for $\tilde{V}$ as in \eqref{cm}, our results in \eqref{Mjklterms} and \eqref{mixed_MSE} give, for the second term of $K_1(\boldsymbol{\eta_0})$, that 
\begin{align}
\label{boundK1}
\nonumber & \frac{1}{2}\left\lbrace\vphantom{(\left(\sup_{\theta:|\theta-\theta_0|\leq\epsilon}\left|l^{(3)}(\theta;\boldsymbol{X})\right|\right)^2}\sum_{k=1}^{2}\sum_{l=1}^{2}\left|\tilde{V}_{lk}\right|\sum_{j=1}^{2}\sum_{i=1}^{2}\sqrt{\EE\left[Q_j^2Q_i^2\right]}\sqrt{\EE\left[\left(nM_{kji}(\boldsymbol{X})\right)^2\middle|\left|Q_{(m)}\right| < \epsilon \right]}\vphantom{(\left(\sup_{\theta:|\theta-\theta_0|\leq\epsilon}\left|l^{(3)}(\theta;\boldsymbol{X})\right|\right)^2}\right\rbrace\\
\nonumber & < \frac{n^2}{\sqrt{2\alpha}(n-5)(n-9)}\left\lbrace\vphantom{(\left(\sup_{\theta:|\theta-\theta_0|\leq\epsilon}\left|l^{(3)}(\theta;\boldsymbol{X})\right|\right)^2}8\left(\sqrt{\eta_1}+|\eta_2|+\eta_1\right)\sqrt{12+\frac{516}{n}+\frac{945}{n^2}}\left(1+\frac{3\left(|\eta_2|+\frac{\eta_1}{2}\right)^2}{\eta_1}\right)\right.\\
\nonumber & \;\left. + 2\sqrt{12(\eta_1+\eta_2^2)^2+\frac{12\eta_2^2(43\eta_2^2+63\eta_1)}{n}+\frac{945\eta_2^4}{n^2}}\right.\\
\nonumber &\left.\qquad\qquad\times\left(\frac{4\left|\eta_2^3\right|}{\eta_1^3} + \frac{(2|\eta_2|+\eta_1)\left(3|\eta_2|+2+2\sqrt{\eta_1}\right)}{\eta_1^2}+\frac{1}{\sqrt{\eta_1}}+1\right)\right.\\
& \;\left. + \frac{4}{\eta_1}\sqrt{4\left(\eta_1+3\eta_2^2\right)+\frac{6}{n}\left(21\eta_1+86\eta_2^2\right)+945\eta_2^2}\left((\eta_1 + |\eta_2|)(\eta_1+3|\eta_2|+2\sqrt{\eta_1})+\eta_1\right)\vphantom{(\left(\sup_{\theta:|\theta-\theta_0|\leq\epsilon}\left|l^{(3)}(\theta;\boldsymbol{X})\right|\right)^2}\right\rbrace,
\end{align}
which is an upper bound for $K_1(\boldsymbol{\eta_0})$. We now proceed to find an upper bound on $K_2(\boldsymbol{\eta_0})$, which is a sum of two quantities as \eqref{remainder_term_R2} shows, involving the calculation of variances of $\xi_{ij}$ as defined in \eqref{cm}. For the first quantity, using \eqref{derivatives_multi_normal} and \eqref{multi_normal_FISHER}, after straightforward calculation of moments (up to fourth order) of $X_1$ and with $\alpha$ and $\beta$ as in the corollary, we get that
\begin{align}
\label{firstTermMultiNormal}
\nonumber & \frac{1}{4}\sum_{j=1}^{2}\left[{\rm Var}\left[\left(\sum_{k=1}^{2}\tilde{V}_{jk}\frac{\partial}{\partial\eta_{k}}\log f(X_1|\boldsymbol{\eta_0})\right)^2\right]\right]^{\frac{1}{2}}\\
\nonumber & = \frac{1}{4}\left\lbrace\sqrt{{\rm Var}\left[\left(X_1^2-\frac{1}{2\eta_1} - \frac{\eta_2^2}{4\eta_1^2}\right)^2\right]}\left(\tilde{V}^2_{11} + \tilde{V}^2_{12}\right) + \sqrt{{\rm Var}\left[\left(X_1-\frac{\eta_2}{2\eta_1}\right)^2\right]}\left(\tilde{V}^2_{22} + \tilde{V}^2_{12}\right)\right\rbrace\\
%\nonumber & = \frac{1}{2\alpha\eta_1^2}\sqrt{\frac{7}{2} + \frac{\eta_2^4}{2\eta_1^2} + \frac{7\eta_2^2}{\eta_1}}\left(\eta_1^3(1+\sqrt{\eta_1})^2 + \eta_1^2\eta_2^2\right) + \frac{1}{2\sqrt{2}\alpha\eta_1}\left(\eta_1^2\eta_2^2 + \left(\eta_1+\eta_1^{\frac{3}{2}}+\eta_2^2\right)^2\right)\\
&= \frac{1}{2\alpha}\left\lbrace \alpha\sqrt{\frac{7}{2} + \frac{\eta_2^4}{2\eta_1^2} + \frac{7\eta_2^2}{\eta_1}}+\frac{1}{\sqrt{2}\eta_1}\left(\eta_1^2\eta_2^2 + \beta^2\right)\right\rbrace.
\end{align}
%\begin{align}
%\nonumber & {\rm Var}\left[\left(\frac{X_1-\mu}{\sigma}\right)^4\right] = \EE\left[\left(\frac{X_1-\mu}{\sigma}\right)^8\right] - \left[\EE\left[\left(\frac{X_1-\mu}{\sigma}\right)^4\right]\right]^2 = 96\\
%\nonumber & {\rm Var}\left[\left(\frac{X_1-\mu}{\sigma}\right)^2\right] = \EE\left[\left(\frac{X_1-\mu}{\sigma}\right)^4\right] - \left[\EE\left[\left(\frac{X_1-\mu}{\sigma}\right)^2\right]\right]^2 =2\\
%\nonumber & {\rm Cov}\left(\left(\frac{X_1-\mu}{\sigma}\right)^4,\left(\frac{X_1-\mu}{\sigma}\right)^2\right) = \EE\left[\left(\frac{X_1-\mu}{\sigma}\right)^6\right]- \EE\left[\left(\frac{X_1-\mu}{\sigma}\right)^4\right]\EE\left[\left(\frac{X_1-\mu}{\sigma}\right)^2\right] =12.
%\end{align}
%Applying the above results to \eqref{firstTermMultiNormalmidre},
%\begin{equation}
%\label{firstTermMultiNormal}
%\frac{\|h\|_2}{4\sqrt{n}}\sum_{j=1}^{d}\left[{\rm Var}\left(\left(\sum_{k=1}^{d}\left[\left[I(\boldsymbol{\theta_0})\right]^{-\frac{1}{2}}\right]_{jk}\frac{\partial}{\partial\theta_{k}}\log f(X_1|\boldsymbol{\theta_0})\right)^2\right)\right]^{\frac{1}{2}} = \frac{\|h\|_2}{4\sqrt{n}}\left(\sqrt{2} + \sqrt{14}\right).
%\end{equation}
For the second quantity in $K_2(\boldsymbol{\eta_0})$, simple calculation of moments leads to
\begin{align}
\label{secondTermMultiNormal}
\nonumber & \frac{1}{2}\left[\vphantom{(\left(\sup_{\theta:|\theta-\theta_0|\leq\epsilon}\left|l^{(3)}(\theta;\boldsymbol{X})\right|\right)^2}{\rm Var}\left[\vphantom{(\left(\sup_{\theta:|\theta-\theta_0|\leq\epsilon}\left|l^{(3)}(\theta;\boldsymbol{X})\right|\right)^2}\sum_{q=1}^{2}\sum_{v=1}^{2}\tilde{V}_{2q}\frac{\partial}{\partial\eta_{q}}\log f(X_1|\boldsymbol{\eta_0})\tilde{V}_{1v}\frac{\partial}{\partial\eta_{v}}\log f(X_1|\boldsymbol{\eta_0})\vphantom{(\left(\sup_{\theta:|\theta-\theta_0|\leq\epsilon}\left|l^{(3)}(\theta;\boldsymbol{X})\right|\right)^2}\right]\vphantom{(\left(\sup_{\theta:|\theta-\theta_0|\leq\epsilon}\left|l^{(3)}(\theta;\boldsymbol{X})\right|\right)^2}\right]^{\frac{1}{2}}\\
\nonumber & \leq \frac{1}{2}\left[\EE\left[\tilde{V}_{11}\tilde{V}_{21}\left(\frac{\partial}{\partial\eta_1}\log f(X_1|\boldsymbol{\eta_0})\right)^2 + \tilde{V}_{22}\tilde{V}_{12}\left(\frac{\partial}{\partial\eta_2}\log f(X_1|\boldsymbol{\eta_0})\right)^2\right.\right.\\
\nonumber & \left.\left.\qquad\qquad\qquad + \frac{\partial}{\partial\eta_1}\log f(X_1|\boldsymbol{\eta_0})\frac{\partial}{\partial\eta_2}\log f(X_1|\boldsymbol{\eta_0})\left(\tilde{V}_{21}^2 + \tilde{V}_{22}\tilde{V}_{11}\right)\right]^2\right]^{\frac{1}{2}}\\
\nonumber & \leq \frac{\sqrt{3}}{2}\left[\tilde{V}^2_{11}\tilde{V}^2_{21}\EE\left[\left(\frac{\partial}{\partial\eta_1}\log f(X_1|\boldsymbol{\eta_0})\right)^4\right] + \tilde{V}_{22}^2\tilde{V}^2_{12}\EE\left[\left(\frac{\partial}{\partial\eta_2}\log f(X_1|\boldsymbol{\eta_0})\right)^4\right]\right.\\
\nonumber & \left.\qquad\qquad\quad + \left(\tilde{V}_{21}^2 + \tilde{V}_{22}\tilde{V}_{11}\right)^2\EE\left[\left(\frac{\partial}{\partial\eta_1}\log f(X_1|\boldsymbol{\eta_0})\frac{\partial}{\partial\eta_2}\log f(X_1|\boldsymbol{\eta_0})\right)^2\right]\right]^{\frac{1}{2}}\\
%\nonumber & = \frac{\sqrt{3}}{2}\left[\frac{3}{\alpha^2}\eta_2^2\left(\alpha-\eta_2^2\right)\left(5+\frac{\eta_2^2}{\eta_1}\right)^2 + \frac{3}{\alpha^2}\eta_2^2\beta^2 + \frac{1}{\alpha^2}\left(2\sqrt{\eta_1}\eta_2^2+\alpha\right)^2\left(5+\frac{3\eta_2^2}{\eta_1}\right)\right]^{\frac{1}{2}}\\
& = \frac{3\eta_2}{2\alpha}\left[\left(\alpha-\eta_2^2\right)\left(5+\frac{\eta_2^2}{\eta_1}\right)^2 + \beta^2+ \left(2\sqrt{\eta_1}\eta_2+\frac{\alpha}{\eta_2}\right)^2\left(\frac{5}{3}+\frac{\eta_2^2}{\eta_1}\right)\right]^{\frac{1}{2}}.
\end{align}
For an upper bound on $K_3(\boldsymbol{\eta_0})$ as in \eqref{remainder_term_R3}, we use that $X_1'$ is an independent copy of $X_1$ and also
\begin{align}
\nonumber & \EE\left[\left|\frac{\partial}{\partial \eta_1}\log f(X_1|\boldsymbol{\eta_0})\right|^3\right] = \EE\left[\left|-X_1^2 + \frac{1}{2\eta_1}+\frac{\eta_2^2}{4\eta_1^2}\right|^3\right]\leq \frac{18}{\eta_1^3}\left(1+\frac{\eta_2^3}{2\eta_1^{\frac{3}{2}}\sqrt{\pi}}\right)\\
\nonumber & \EE\left[\left|\frac{\partial}{\partial \eta_2}\log f(X_1|\boldsymbol{\eta_0})\right|^3\right] = \frac{1}{\sqrt{\pi}\eta_1^{\frac{3}{2}}}.
\end{align}
Then, the triangle inequality and \eqref{interestinginequality3} yield
\begin{align}
\label{thirdTermMultiNormal}
\nonumber & \frac{1}{12}\EE\left[\sum_{i=1}^{2}\left|\sum_{l=1}^{2}\tilde{V}_{il}\left(\frac{\partial}{\partial \eta_l}\log f(X_1'|\boldsymbol{\eta_0}) - \frac{\partial}{\partial \eta_l}\log f(X_1|\boldsymbol{\eta_0})\right)\right|\right]^3\\
\nonumber & \leq \frac{32}{3}\left\lbrace\EE\left[\left|\frac{\partial}{\partial \eta_1}\log f(X_1|\boldsymbol{\eta_0})\right|^3\right]\left(\left|\tilde{V}_{11}\right|^3 + \left|\tilde{V}_{21}\right|^3\right)\right.\\
\nonumber & \left.\qquad\qquad\qquad + \EE\left[\left|\frac{\partial}{\partial \eta_2}\log f(X_1|\boldsymbol{\eta_0})\right|^3\right]\left(\left|\tilde{V}_{12}\right|^3 + \left|\tilde{V}_{22}\right|^3\right)\right\rbrace\\
& = \frac{64\sqrt{2}}{3\alpha^{\frac{3}{2}}}\left\lbrace\vphantom{(\left(\sup_{\theta:|\theta-\theta_0|\leq\epsilon}\left|l^{(3)}(\theta;\boldsymbol{X})\right|\right)^2}18\left(1+\frac{\eta_2^3}{2\eta_1^{\frac{3}{2}}\sqrt{\pi}}\right)\left(\eta_1^{\frac{3}{2}}\left(1+\sqrt{\eta_1}\right)^3 + \left|\eta_2\right|^3\right)+ \frac{\eta_1^3|\eta_2|^3 + \beta^3}{\sqrt{\pi}\eta_1^{\frac{3}{2}}}\vphantom{(\left(\sup_{\theta:|\theta-\theta_0|\leq\epsilon}\left|l^{(3)}(\theta;\boldsymbol{X})\right|\right)^2}\right\rbrace.
\end{align}
For the last term of \eqref{final_bound_regression}, we obtain that
\begin{align}
\label{fourthTermMultiNormal}
\nonumber \frac{2\|h\|}{\epsilon^2}\EE\left[\sum_{j=1}^{2}Q_j^2\right] &= \frac{2\|h\|}{\epsilon^2(n-3)(n-5)}\left((2n+15)\eta_1^2+2n\left(\eta_2^2+\eta_1\right)+15\eta_2^2\right)\\
& = \frac{8\|h\|}{\eta_1^2(n-3)(n-5)}\left((2n+15)\left(\eta_1^2 + \eta_2^2\right)+2n\eta_1\right),
\end{align}
where for the second equality we used that $\epsilon = \frac{\eta_1}{2}$, with our choice explained in the paragraph after \eqref{Mjklterms}. Using the results in \eqref{Tkjzero}, \eqref{boundK1}, \eqref{firstTermMultiNormal}, \eqref{secondTermMultiNormal}, \eqref{thirdTermMultiNormal} and \eqref{fourthTermMultiNormal} we get the result of the corollary. $\;\;\;\qquad\qquad\qquad\qquad\qquad\qquad\qquad\qquad\qquad\qquad\qquad\qquad\qquad\qquad\;\;\;\qquad\quad\;\;\;\;\;\;\;\;\blacksquare$
\vspace{0.1in}
\\
{\raggedright \textit{Proof}{\rm of \textbf{Corollary \ref{Cor_multi_Beta}}.}}\\
Part \textbf{a)}. The probability density function is
\begin{equation}
\nonumber f(x|\boldsymbol{\theta}) = \frac{\Gamma(\alpha+\beta)}{\Gamma(\alpha)\Gamma(\beta)}x^{\alpha-1}(1-x)^{\beta-1},
\end{equation}
with $\alpha, \beta >0$ and $x \in [0,1]$. Hence, for $j,k\in \mathbb{Z^+}$
\begin{align}
\label{multi_Beta_likelihood}
\nonumber & \frac{\partial}{\partial \alpha} \log f(x|\boldsymbol{\theta}) = \Psi(\alpha + \beta) - \Psi(\alpha) + \log(x),\\
\nonumber & \frac{\partial}{\partial \beta} \log f(x|\boldsymbol{\theta}) = \Psi(\alpha + \beta) - \Psi(\beta) + \log(1-x)\\
\nonumber & \frac{\partial^{j+1}}{\partial\alpha^{j+1}} \log f(x|\boldsymbol{\theta}) = \Psi_{j}(\alpha + \beta) - \Psi_{j}(\alpha),\\
\nonumber & \frac{\partial^{j+1}}{\partial\beta^{j+1}} \log f(x|\boldsymbol{\theta}) = \Psi_{j}(\alpha + \beta) - \Psi_{j}(\beta)\\
& \frac{\partial^{k+j}}{\partial \alpha^k \partial \beta^j}\log f(x|\boldsymbol{\theta}) = \Psi_{k+j-1}(\alpha + \beta).
\end{align}
From \eqref{multi_Beta_likelihood}, we see that we are under the scenario \textbf{(2)} of Remark \ref{remark_implicit} and $U_1$ will be calculated using \eqref{special_case}. The expected Fisher Information matrix is
\begin{equation}
\nonumber I(\boldsymbol{\theta_0}) = \begin{pmatrix}
\Psi_1(\alpha) - \Psi_1(\alpha + \beta) & -\Psi_1(\alpha + \beta)\\
-\Psi_1(\alpha + \beta) & \Psi_1(\beta) - \Psi_1(\alpha + \beta)
\end{pmatrix}.
\end{equation}
Simple calculations show that the inverse of $I(\boldsymbol{\theta_0})$ is
\begin{equation}
\nonumber \left[I(\boldsymbol{\theta_0})\right]^{-1} = \frac{1}{\delta_I}\begin{pmatrix}
C_1(\beta,\alpha) & \Psi_1(\alpha + \beta)\\
\Psi_1(\alpha + \beta) & C_1(\alpha,\beta)
\end{pmatrix}.
\end{equation}
Therefore,
\begin{equation}
%\label{root_Fisher_Beta}
\nonumber \left[I(\boldsymbol{\theta_0})\right]^{-2} = \frac{1}{\delta^2_I}\begin{pmatrix}
C^2_1(\beta,\alpha) + \Psi^2_1(\alpha+\beta) & \Psi_1(\alpha + \beta)(C_1(\alpha,\beta) + C_1(\beta,\alpha))\\
\Psi_1(\alpha + \beta)(C_1(\alpha,\beta) + C_1(\beta,\alpha)) & C^2_1(\alpha,\beta) + \Psi^2_1(\alpha+\beta)
\end{pmatrix}.
\end{equation}
For $k, q =1,2$, we now proceed to calculate the quantities
\begin{equation}
\nonumber \mathop{\sup_{\boldsymbol{\theta}:\left|\theta_j - \theta_{0,j}\right|<\epsilon}}_{\forall j \in \left\lbrace 1,2\right\rbrace}\left\lbrace\left[I^{-2}(\boldsymbol{\theta})\right]_{kq}\right\rbrace.
\end{equation}
Firstly, the fact that $\delta_I(\alpha, \beta)$ as in \eqref{delta} is a positive, decreasing function of $\alpha$ and $\beta$, means that
\begin{equation}
\label{bounddelta}
\mathop{\sup_{\boldsymbol{\theta}:\left|\theta_j - \theta_{0,j}\right|<\epsilon}}_{\forall j \in \left\lbrace 1,2\right\rbrace}\left\lbrace\frac{1}{[\delta_I(\theta_1,\theta_2)]^2}\right\rbrace = \frac{1}{[\delta_I(\alpha+\epsilon,\beta+\epsilon)]^2}.
\end{equation}
In regards to $C_1^2(\theta_1,\theta_2)$ as in \eqref{delta}, we have that using a first-order Taylor expansion and for $\tilde{\theta}$ between $\theta_1$ and $\theta_1+\theta_2$,
\begin{align}
\label{boundC1}
\nonumber \mathop{\sup_{\boldsymbol{\theta}:\left|\theta_j - \theta_{0,j}\right|<\epsilon}}_{\forall j \in \left\lbrace 1,2\right\rbrace}\left\lbrace C_1^2(\theta_1,\theta_2)\right\rbrace &= \mathop{\sup_{\boldsymbol{\theta}:\left|\theta_j - \theta_{0,j}\right|<\epsilon}}_{\forall j \in \left\lbrace 1,2\right\rbrace}\left\lbrace\theta_2^2\Psi^2_2(\tilde{\theta})\right\rbrace\\
& = (\beta+\epsilon)^2\mathop{\sup_{\boldsymbol{\theta}:\left|\theta_j - \theta_{0,j}\right|<\epsilon}}_{\forall j \in \left\lbrace 1,2\right\rbrace}\left\lbrace\Psi_2^2(\theta_1)\right\rbrace = (\beta+\epsilon)^2\Psi_2^2(\alpha - \epsilon),
\end{align}
since $\Psi_2^2(x)$ is a decreasing function of $x$; see the definition of $\Psi_2(\cdot)$ in \eqref{psi_m}. In the same way, we can find an upper bound for $C_1^2(\theta_2,\theta_1)$. With regards to the quantity $C_1(\theta_1 + \theta_2) + C_1(\theta_2 + \theta_1)$, we have that a similar first-order Taylor expansion as in \eqref{boundC1} leads to
\begin{equation}
\label{boundC1C1mid}
C_1(\theta_1 + \theta_2) + C_1(\theta_2 + \theta_1) = -\theta_2\Psi_2\left(\tilde{\theta}\right) - \theta_1\Psi_2\left(\tilde{\tilde{\theta}}\right),
\end{equation}
where $\tilde{\theta}$ is between $\theta_1$ and $\theta_1 + \theta_2$, while $\tilde{\tilde{\theta}}$ is between $\theta_2$ and $\theta_1 + \theta_2$. It is important to highlight that $\Psi_2(x)$, as defined in \eqref{psi_m} is a negative and increasing function of $x$. Continuing from \eqref{boundC1C1mid},
\begin{equation}
\label{boundC1C1}
\mathop{\sup_{\boldsymbol{\theta}:\left|\theta_j - \theta_{0,j}\right|<\epsilon}}_{\forall j \in \left\lbrace 1,2\right\rbrace}\left\lbrace C_1(\theta_1,\theta_2) + C_1(\theta_2,\theta_1)\right\rbrace = -\left[\left(\beta+\epsilon\right)\Psi_2(\alpha-\epsilon) + \left(\alpha+\epsilon\right)\Psi_2(\beta-\epsilon)\right].
\end{equation}
Using the results in \eqref{bounddelta}, \eqref{boundC1}, as well as the fact that $\Psi_1(x)$ defined in \eqref{psi_m} is a positive, decreasing function of $x$, we have that
\begin{align}
\label{supremumsJ}
\nonumber & \mathop{\sup_{\boldsymbol{\theta}:\left|\theta_j - \theta_{0,j}\right|<\epsilon}}_{\forall j \in \left\lbrace 1,2\right\rbrace}\left\lbrace\left[I^{-2}(\boldsymbol{\theta})\right]_{11}\right\rbrace = \frac{(\alpha+\epsilon)^2\Psi_2^2(\beta - \epsilon) + \Psi_1^2(\alpha + \beta - 2\epsilon)}{[\delta_I(\alpha+\epsilon,\beta+\epsilon)]^2}\\
\nonumber & \mathop{\sup_{\boldsymbol{\theta}:\left|\theta_j - \theta_{0,j}\right|<\epsilon}}_{\forall j \in \left\lbrace 1,2\right\rbrace}\left\lbrace\left[I^{-2}(\boldsymbol{\theta})\right]_{12}\right\rbrace = \mathop{\sup_{\boldsymbol{\theta}:\left|\theta_j - \theta_{0,j}\right|<\epsilon}}_{\forall j \in \left\lbrace 1,2\right\rbrace}\left\lbrace\left[I^{-2}(\boldsymbol{\theta})\right]_{21}\right\rbrace\\
\nonumber & \qquad\qquad\qquad\qquad\qquad\quad = -\frac{\left[\left(\beta+\epsilon\right)\Psi_2(\alpha-\epsilon) + \left(\alpha+\epsilon\right)\Psi_2(\beta-\epsilon)\right]}{[\delta_I(\alpha+\epsilon,\beta+\epsilon)]^2}\\
& \mathop{\sup_{\boldsymbol{\theta}:\left|\theta_j - \theta_{0,j}\right|<\epsilon}}_{\forall j \in \left\lbrace 1,2\right\rbrace}\left\lbrace\left[I^{-2}(\boldsymbol{\theta})\right]_{22}\right\rbrace = \frac{(\beta+\epsilon)^2\Psi_2^2(\alpha - \epsilon) + \Psi_1^2(\alpha + \beta - 2\epsilon)}{[\delta_I(\alpha+\epsilon,\beta+\epsilon)]^2}.
\end{align}
To derive the expression for $U_1$ as in \eqref{special_case} in this special case, we need to calculate the quantities $\EE\left[\frac{\partial}{\partial\theta_k}\ell(\boldsymbol{\theta_0};\boldsymbol{X})\frac{\partial}{\partial\theta_q}\ell(\boldsymbol{\theta_0};\boldsymbol{X})\right]$, for $k,q =1,2$. Using \eqref{multi_Beta_likelihood}, we have that
\begin{align}
\label{expectationscorefunction}
\nonumber & \EE\left[\left(\frac{\partial}{\partial\alpha}\ell(\boldsymbol{\theta_0};\boldsymbol{X})\right)^2\right] = \EE\left[\left(n(\Psi(\alpha + \beta) - \Psi(\alpha)) + \sum_{i=1}^n\log(X_i)\right)^2\right]\\
\nonumber & \qquad\qquad\qquad\qquad\quad = {\rm Var}\left[\sum_{i=1}^n\log(X_i)\right] = nC_1(\alpha,\beta)\\
\nonumber & \EE\left[\left(\frac{\partial}{\partial\beta}\ell(\boldsymbol{\theta_0};\boldsymbol{X})\right)^2\right] = \EE\left[\left(n(\Psi(\alpha + \beta) - \Psi(\beta)) + \sum_{i=1}^n\log(1-X_i)\right)^2\right]\\
\nonumber & \qquad\qquad\qquad\qquad\quad = {\rm Var}\left[\sum_{i=1}^n\log(1-X_i)\right] = nC_1(\beta,\alpha)\\
\nonumber & \EE\left[\frac{\partial}{\partial\alpha}\ell(\boldsymbol{\theta_0};\boldsymbol{X})\frac{\partial}{\partial\beta}\ell(\boldsymbol{\theta_0};\boldsymbol{X})\right] = n(\Psi(\alpha+\beta) - \Psi(\beta))\EE\left[\sum_{i=1}^n\log(X_i)\right] + \EE\left[\sum_{i=1}^n\sum_{j=1}^n\log(X_i)\log(1-X_j)\right]\\
\nonumber & \qquad\qquad\qquad\qquad\quad\quad\qquad\; = n^2(\Psi(\alpha+\beta) - \Psi(\beta))(\Psi(\alpha) - \Psi(\alpha+\beta))\\
\nonumber & \qquad\qquad\qquad\qquad\qquad\qquad\quad + n\left((\Psi(\alpha) - \Psi(\alpha+\beta))(\Psi(\beta) - \Psi(\alpha+\beta)) - \Psi_1(\alpha+\beta)\right)\\
\nonumber & \qquad\qquad\qquad\qquad\qquad\qquad\quad + n(n-1)(\Psi(\alpha) - \Psi(\alpha+\beta))(\Psi(\beta) - \Psi(\alpha+\beta))\\
& \qquad\qquad\qquad\qquad\quad\quad\qquad\; = -n\Psi_1(\alpha+\beta).
\end{align}
Applying the results of \eqref{supremumsJ} and \eqref{expectationscorefunction} to \eqref{special_case}, we conclude that
\begin{align}
\nonumber \EE\left[\sum_{j=1}^2Q_j^2\right] \leq & \frac{1}{n[\delta_I(\alpha+\epsilon,\beta+\epsilon)]^2}\left\lbrace\vphantom{(\left(\sup_{\theta:|\theta-\theta_0|\leq\epsilon}\left|l^{(3)}(\theta;\boldsymbol{X})\right|\right)^2}C_1(\alpha,\beta)\left[(\alpha+\epsilon)^2\Psi_2^2(\beta - \epsilon) + \Psi_1^2(\alpha + \beta - 2\epsilon)\right]\right.\\
\nonumber & \left. + C_1(\beta,\alpha)\left[(\beta+\epsilon)^2\Psi_2^2(\alpha - \epsilon) + \Psi_1^2(\alpha + \beta - 2\epsilon)\right]\right.\\
\nonumber & \left. + 2\Psi_1(\alpha+\beta)\left[\left(\beta+\epsilon\right)\Psi_2(\alpha-\epsilon) + \left(\alpha+\epsilon\right)\Psi_2(\beta-\epsilon)\right]\vphantom{(\left(\sup_{\theta:|\theta-\theta_0|\leq\epsilon}\left|l^{(3)}(\theta;\boldsymbol{X})\right|\right)^2}\right\rbrace,
\end{align}
which completes the proof.$\qquad\qquad\qquad\qquad\;\;\;\;\;\;\;\qquad\qquad\qquad\qquad\qquad\qquad\qquad\qquad\qquad\quad\blacksquare$

\section*{Acknowledgements}
This research occurred whilst Andreas Anastasiou was at the University of Oxford, supported by a Teaching Assistantship Bursary from the Department of Statistics, University of Oxford, and the Engineering and Physical Sciences Research Council (EPSRC) grant EP/K503113/1. The author would like to thank Gesine Reinert, Tobias Kley, and Christophe Ley for insightful comments.

\bibliographystyle{Chicago}
\bibliography{Revised_EJS_paper_general}

\begin{thebibliography}{}

\bibitem[\protect\citeauthoryear{Anastasiou}{Anastasiou}{2017}]{Anastasiou_m_dependence}
Anastasiou, A. (2017).
\newblock Bounds for the normal approximation of the maximum likelihood
  estimator from $m-$dependent random variables.
\newblock {\em Statistics \& Probability Letters\/}~{\em {\bf{129}}}, 171--181.

\bibitem[\protect\citeauthoryear{Anastasiou and Ley}{Anastasiou and
  Ley}{2017}]{Anastasiou_Ley}
Anastasiou, A. and C.~Ley (2017).
\newblock Bounds for the asymptotic normality of the maximum likelihood
  estimator using the {D}elta method.
\newblock {\em ALEA, Lat. Am. J. Probab. Math. Stat.\/}~{\em {\bf{14}}},
  153--171.

\bibitem[\protect\citeauthoryear{Anastasiou and Reinert}{Anastasiou and
  Reinert}{2017}]{Anastasiou_Reinert}
Anastasiou, A. and G.~Reinert (2017).
\newblock Bounds for the normal approximation of the maximum likelihood
  estimator.
\newblock {\em Bernoulli\/}~{\em {\bf 23}}, 191--218.

\bibitem[\protect\citeauthoryear{Berk}{Berk}{1972}]{Berk}
Berk, R.~H. (1972).
\newblock Consistency and asymptotic normality of {MLE'}s for exponential
  models.
\newblock {\em The Annals of Mathematical Statistics\/}~{\em {\bf 43}},
  193--204.

\bibitem[\protect\citeauthoryear{Billingsley}{Billingsley}{1961}]{Billingsley}
Billingsley, P. (1961).
\newblock Statistical {M}ethods in {M}arkov {C}hains.
\newblock {\em The Annals of Mathematical Statistics\/}~{\em {\bf 32}, No.1},
  12--40.

\bibitem[\protect\citeauthoryear{Casella and Berger}{Casella and
  Berger}{2002}]{Casella}
Casella, G. and R.~L. Berger (2002).
\newblock {\em Statistical {I}nference\/} (Second ed.).
\newblock Brooks/Cole, Cengage Learning, Duxbury, Pacific Grove.

\bibitem[\protect\citeauthoryear{Davison}{Davison}{2008}]{Davison}
Davison, A.~C. (2008).
\newblock {\em Statistical {M}odels\/} (First ed.).
\newblock Cambridge Series in Statistical and Probabilistic Mathematics.
  Cambridge University Press.

\bibitem[\protect\citeauthoryear{Fahrmeir and Kaufmann}{Fahrmeir and
  Kaufmann}{1985}]{Fahrmeir}
Fahrmeir, L. and H.~Kaufmann (1985).
\newblock Consistency and asymptotic normality of the maximum likelihood
  estimator in generalized linear models.
\newblock {\em The Annals of Statistics\/}~{\em {\bf 13}, No.1}, 342--368.

\bibitem[\protect\citeauthoryear{Hoadley}{Hoadley}{1971}]{Hoadley}
Hoadley, B. (1971).
\newblock Asymptotic {P}roperties of {M}aximum {L}ikelihood {E}stimators for
  the {I}ndependent {N}ot {I}dentically {D}istributed {C}ase.
\newblock {\em The Annals of Mathematical Statistics\/}~{\em {\bf 42}, No.6},
  1977--1991.

\bibitem[\protect\citeauthoryear{Kiefer}{Kiefer}{1968}]{Kiefer}
Kiefer, J.~C. (1968).
\newblock Statistical inference.
\newblock In {\em The future of statistics. Proceedings of a {C}onference on
  the {F}uture of {S}tatistics held at the {U}niversity of {W}isconsin,
  {M}adison, {W}isconsin, {J}une 1967}, pp.\  139--142. Academic Press, New
  York-London.

\bibitem[\protect\citeauthoryear{Koroljuk and Borovskich}{Koroljuk and
  Borovskich}{1994}]{Koroljuk}
Koroljuk, V.~S. and Y.~V. Borovskich (1994).
\newblock {\em Theory of U-statistics}.
\newblock Mathematics and its Applications {\textbf{273}}. Kluwer Academic
  Publishers Group, Dordrecht. Translated from the 1989 Russian original by
  P.V. Malyshev and D.V. Malyshev and revised by the authors.

\bibitem[\protect\citeauthoryear{Lauritzen}{Lauritzen}{1988}]{Steffen}
Lauritzen, S. (1988).
\newblock {\em Extremal {F}amilies and {S}ystems of {S}ufficient {S}tatistics}.
\newblock Lecture Notes in Statistics, No.49. Springer-Verlag,
  Berlin-Heidelberg-New York.

\bibitem[\protect\citeauthoryear{Lauritzen}{Lauritzen}{1996}]{Lauritzen_graphical}
Lauritzen, S. (1996).
\newblock {\em Graphical {M}odels}.
\newblock Oxford: Clarendon Press.

\bibitem[\protect\citeauthoryear{M\"akel\"ainen, Schmidt, and
  Styan}{M\"akel\"ainen et~al.}{1981}]{Makelainen}
M\"akel\"ainen, T., K.~Schmidt, and G.~P.~H. Styan (1981).
\newblock On the existence and uniqueness of the maximum likelihood estimate of
  a vector-valued parameter in fixed-size samples.
\newblock {\em The Annals of Statistics\/}~{\em {\bf 9}, No.4}, 758--767.

\bibitem[\protect\citeauthoryear{Massam, Li, and Gao}{Massam
  et~al.}{2018}]{Massam}
Massam, H., Q.~Li, and X.~Gao (2018).
\newblock Bayesian precision and covariance matrix estimation for graphical
  {G}aussian models with edge and vertex symmetries.
\newblock {\em Biometrika, asx084, https://doi.org/10.1093/biomet/asx084\/},
  1--18.

\bibitem[\protect\citeauthoryear{Pinelis}{Pinelis}{2017}]{Pinelis}
Pinelis, I. (2017).
\newblock Optimal-order uniform and nonuniform bounds on the rate of
  convergence to normality for maximum likelihood estimators.
\newblock {\em Electronic Journal of Statistics\/}~{\em {\bf{11}}}, 1160--1179.

\bibitem[\protect\citeauthoryear{Pinelis and Molzon}{Pinelis and
  Molzon}{2016}]{Pinelis_Molzon}
Pinelis, I. and R.~Molzon (2016).
\newblock Optimal-order bounds on the rate of convergence to normality in the
  multivariate delta method.
\newblock {\em Electronic Journal of Statistics\/}~{\em {\bf{10}}}, 1001--1063.

\bibitem[\protect\citeauthoryear{Reinert and R\"{o}llin}{Reinert and
  R\"{o}llin}{2009}]{ReinertRollin}
Reinert, G. and A.~R\"{o}llin (2009).
\newblock Multivariate normal approximation with {S}tein's method of
  exchangeable pairs under a general linearity condition.
\newblock {\em The Annals of Probability\/}~{\em {\bf 37}, No.6}, 2150--2173.

\bibitem[\protect\citeauthoryear{Stein}{Stein}{1972}]{Stein1972}
Stein, C. (1972).
\newblock A bound for the error in the normal approximation to the distribution
  of a sum of dependent random variables.
\newblock In {\em Proceedings of the {S}ixth {B}erkeley {S}ymposium on
  {M}athematical {S}tatistics and {P}robability}, Volume~{\bf 2}, pp.\
  586--602. Berkeley: University of California Press.

\end{thebibliography}

\end{document}